\documentclass[a4paper,11pt]{article}
\usepackage{graphicx}
\usepackage{amsfonts}
\usepackage{amsmath}
\usepackage{amssymb}
\usepackage{indentfirst}
\usepackage{titlesec}

\begin{document}

\renewcommand{\abstractname}{Abstract}

\renewcommand{\figurename}{Fig.}

\providecommand{\ve}{\varepsilon}

\providecommand{\play}{\displaystyle}

\providecommand{\li}{\limits}

\providecommand{\ra}{\rightarrow}

\providecommand{\da}{\downarrow}

\title{\textbf{On stochasticity in nearly-elastic systems}}
\author{Mark Freidlin\thanks{mif@math.umd.edu} , Wenqing Hu\thanks{huwenqing@math.umd.edu}\\Department
of Mathematics, University of Maryland at College Park}
\date{}
\maketitle

\begin{abstract}
Nearly-elastic model systems with one or two degrees of freedom are
considered: the system is undergoing a small loss of energy in each
collision with the "wall". We show that instabilities in this purely
deterministic system lead to stochasticity of its long-time
behavior. Various ways to give a rigorous meaning to the last
statement are considered. All of them, if applicable, lead to the
same stochasticity which is described explicitly. So that the
stochasticity of the long-time behavior is an intrinsic property of
the deterministic systems.
\end{abstract}



\section{Introduction}

Consider the one-dimensional motion of a unit-mass particle in a
smooth potential $F(q)$ in an interval $[a_1,a_2]$ with elastic
reflection at $a_1$ and $a_2$. Let $F'(a_1)>0$ and $F'(a_2)<0$
(Fig.1(a), the potential $F(q)$ is shown there as the bold-face
line): if $q\in [a_1,a_2)$ and $p>0$ or $q\in (a_1,a_2]$ and $p<0$,
the trajectory starting at $(q,p)$ moves according to equation
$\ddot{q}(t)=-F'(q(t))$; the trajectory jumps instantaneously from
$(a_i,p)$ to $(a_i,-p)$, if $i=2$ and $p>0$ or if $i=1$ and $p<0$.
If the initial velocity is large enough, the particle hits both
"walls" $a_1$ and $a_2$ and performs periodic oscillations shown in
Fig.1(b). Let now the walls be not absolutely elastic. It is natural
to assume in certain situations that the loss of energy is
proportional to the speed at the collision point with the wall: if
the particle hits $a_i$ with a speed $v$, it is instantly reflected
with the speed $-v(1-\ve c_i(v))$, $i\in \{1,2\}$. Here $c_i(v)$ are
positive smooth functions, $0<\ve<<1$. The coefficient $(1-\ve c_i
(v))$ is called the coefficient of restitution. Denote by $q^\ve(t)$
the position at time $t$ of the particle performing nearly-elastic
motion, $p^\ve(t)=\dot{q}^\ve(t)$. It is clear that, for each $\ve
>0$ and $t$ large enough, $q^\ve (t)$ will be situated in an
arbitrary small neighborhood of one of the points $a_1, a_2, a_3$
(We assume that the potential $F(q)$ has the form shown in Fig.1(a)
and $a_3$ is the unique local minimum inside $[a_1, a_2]$).

\begin{figure}
\centering
\includegraphics[height=9cm, width=15cm , bb=60 4 537 305]{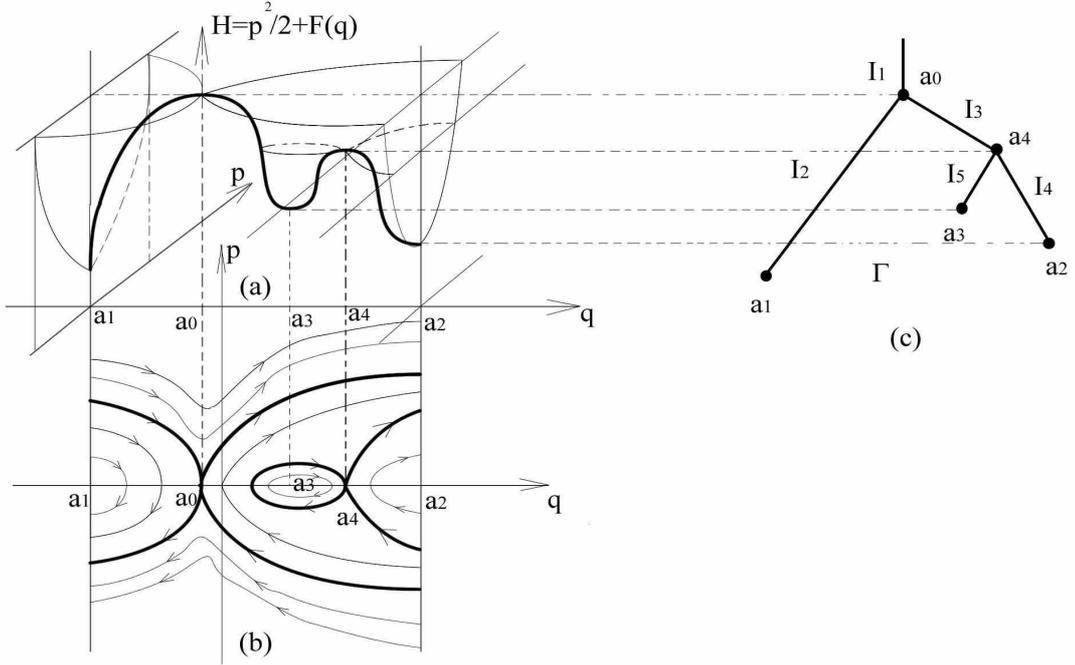}
\caption{The 1-dimensional mechanical model}
\end{figure}

To be specific, assume that the initial speed
$\dot{q}^\ve(0)=p^\ve(0)=v_0$ is large enough : $v_0 > \sqrt{2\max
F(q)}$. Then $q^\ve(t)$, for $\ve$ small enough, hits both points
$a_1$ and $a_2$, but $\lim\limits_{t \rightarrow
\infty}q^\ve(t)=q^\ve_{\infty}$ exists and $q^\ve_{\infty}=a_1$ or
$q^\ve_{\infty}=a_2$. Which of these two points appears as
$\lim\limits_{t \rightarrow \infty}q^\ve(t)$ depends on $\ve$ in a
very sensitive way so that $\lim\li_{\ve \da
0}q^\ve_{\infty}=\lim\li_{\ve \da 0}\lim\li_{t \ra
\infty}q^{\ve}(t)$ does not exist. We show that the final position
$q^\ve_{\infty}$ converges, in a certain sense, as $\ve \da 0$ to a
random variable distributed between the points $a_1$ and $a_2$.
There are various ways to make the last statement rigorous. But they
all lead to the same distribution between $a_1$ and $a_2$, so that
the stochasticity of the system as $\ve << 1$ is an intrinsic
property of the system.

The perturbed system $q^\ve(t)$ for $0 <\ve << 1$ has fast and slow
components. The fast component consists of the motion along the
non-perturbed trajectory. To describe the slow component, consider
the graph $\Gamma$ obtained after identification of points of each
connected component of every level set of the Hamiltonian
$H(q,p)=\play{\frac{p^2}{2}+F(q)}$ (Fig.1(c)). Denote by $Y: \sqcap
\ra \Gamma$ the identification map of the phase space
$\sqcap=\{(q,p)\in \mathbb{R}^2: a_1\leq q \leq a_2\}$ of our system
on $\Gamma$. The slow component of the motion is
$Y(q^\ve(t),\dot{q}^\ve(t))$ (compare with [6], Ch.8). Number the
edges of the graph ($\Gamma= I_1\cup I_2\cup ... \cup I_5$ in the
Fig.1(c)). Then $Y(q,p)=(H(q,p), K(q,p))$, $(q,p)\in \sqcap$, where
$K(q,p)$ is the number of the edge containing $Y(q,p)$ and
$H(q,p)=\play{\frac{p^2}{2}+F(q)}$. The pair $(H,K)$ form a global
coordinate system on $\Gamma$, and the slow component is the motion
$Y_t^\ve=(H(q^\ve(t),\dot{q}^\ve(t)), K(q^\ve(t),\dot{q}^\ve(t)))$
on $\Gamma$.

We prove that, in a certain sense, the rescaled slow motion
$Y^\ve_{t/\ve}$ converges weakly to a stochastic process $Y_t$ on
$\Gamma$. Inside the edges, $Y_t$ is a deterministic motion and the
convergence follows from more or less standard averaging principle
(see, for instance, [1], Ch.10). The stochasticity appears due to a
branching at the interior vertices of $\Gamma$ (compare [2]).

Note that in the case of system shown in Fig.1, the phase trajectory
$(q^\ve(t),\dot{q}^\ve(t))$ with
$\dot{q}^\ve_0=v_0>\sqrt{2\max\li_{a_1\leq q \leq a_2}F(q)}$ never
enters the domain $\{(q,p)\in \sqcap: Y(q,p)\in I_5\}$, so that if
the damping of the system occurs just on the walls $a_1$ and $a_2$,
the stochasticity of the limiting slow motion is concentrated at the
vertex corresponding to the absolute maximum of the potential.
Therefore we will consider in more detail the case when $F(q)$ has
just one local maximum on $[a_1,a_2]$. If the potential has many
wells and the system is losing energy not just at the walls $a_1$,
$a_2$, one should take into account the whole graph $\Gamma$ even if
the initial speed is large.

In the next two sections we consider a model problem where the
limiting slow motion inside the edges is just a motion with constant
speed. This system has some damping not just in the ends of the
interval $[a_1,a_2]$ and models systems with multiwell Hamiltonian.
To give a rigorous meaning to the statement that the slow motion
converges to a stochastic process, we, first, perturb the system
$q^\ve(t)$ by a stochastic perturbation of a small intensity
$\delta<<1$. Then in this already stochastic system
$q^{\ve,\delta}(t)$, we consider the slow component
$Y_t^{\ve,\delta}=Y(q^{\ve,\delta}(t), \dot{q}^{\ve,\delta}(t))$. We
show that $Y_{t/\ve}^{\ve,\delta}$ converges, as first $\ve \da 0$
 and then $\delta \da 0$, in a certain sense to a stochastic process
$Y_t$ on $\Gamma$. This double limit $\lim\li_{\delta \da
0}\lim\li_{\ve \da 0}Y^{\ve,\delta}_t=Y_t$ exists and is the same
for a broad class of stochastic regularizations. It is in this sense
that the slow component of the purely deterministic system
$q^\ve(t)$ should be approximated by the stochastic process $Y_t$ on
$\Gamma$.

\begin{figure}
\centering
\includegraphics[height=6cm, width=10cm , bb=20 2 408 266]{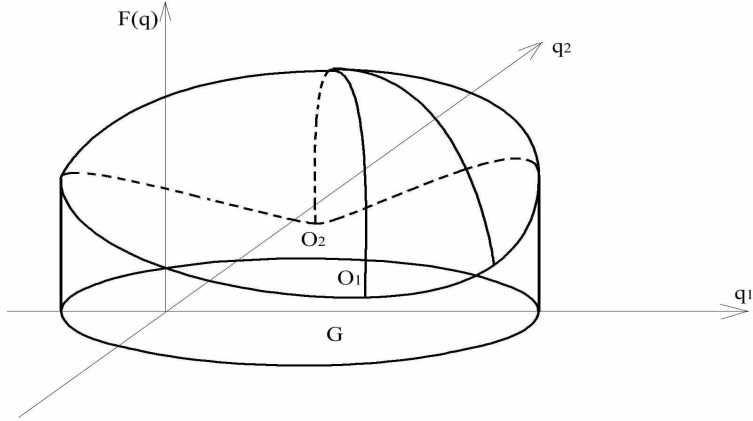}
\caption{A 2-dimensional problem}
\end{figure}

In section 4 we come back to the problem mentioned in the beginning
of this section.

Of course, one can consider similar questions for systems with more
than one degree of freedom. Let, say, $F(q)$, $q\in \mathbb{R}^2$ be
as shown in Fig.2. The particle of a unit mass moves inside a convex
domain $G$ with a smooth boundary $\partial G$ according to the
equation $\ddot{q}_t=-\nabla F(q_t)$ and undergoes an instantaneous
mirror reflection on the boundary. If the collisions with the wall
$\partial G$ are absolutely elastic, the particle will move forever
on the energy surface $H(p,q)=\play{\frac{|p|^2}{2}+F(q)}=H_0$. But
if the particle loses energy in each collision similar to the case
of one degree of freedom, it can, eventually, be found near one of
the local minima of $F(q)$ on $\partial G$ (points $O_1$ and $O_2$
in Fig.2). We consider many-degrees-of-freedom problems for model
potentials in Section 5.

\section{One-degree-of-freedom Model Problems}

Consider a model of the system with several potential wells. Suppose
a particle of unit mass moves freely in an interval $[q_1,q_n]$ with
elastic reflection at the ends of the interval if the initial
velocity is large enough. Let a finite number of points
$q_2,q_3,...,q_{n-1}\in (q_1,q_n)$ and $0<p_2<...<p_{n-1}$ are given
such that if the particle starts with a velocity $p_0$ at a point
$q_0\in [q_1,q_n]$ then it moves freely in
$[q_{i_-(p_0,q_0)},q_{i_+(p_0,q_0)}]$ with instantaneous reflection
in the ends of this interval, where $i_-=i_-(p_0,q_0)$ and
$i_+=i_+(p_0,q_0)$ are defined by the conditions (see Fig.3):
$$q_{i_-}<q_0<q_{i_+} \ \ , \ \ p_0\leq \min(p_{i_-},p_{i_+}) \ \ , \ \ p_0>p_k \ \ for \ \ k\in \{i: q_i\in (q_{i_-},q_{i_+})\}. $$
The energy $\displaystyle{\frac{p_0^2}{2}}$ is preserved in this
system. This is an approximation of the motion in a potential which
has sharp maxima at points $q_i$ and is close to a constant between
$q_i$ and $q_{i+1}$. Just for brevity, we assume that the walls
between the maxima have the same depth.

\begin{figure}
\centering
\includegraphics[height=8cm, width=10cm , bb=1 25 439 305]{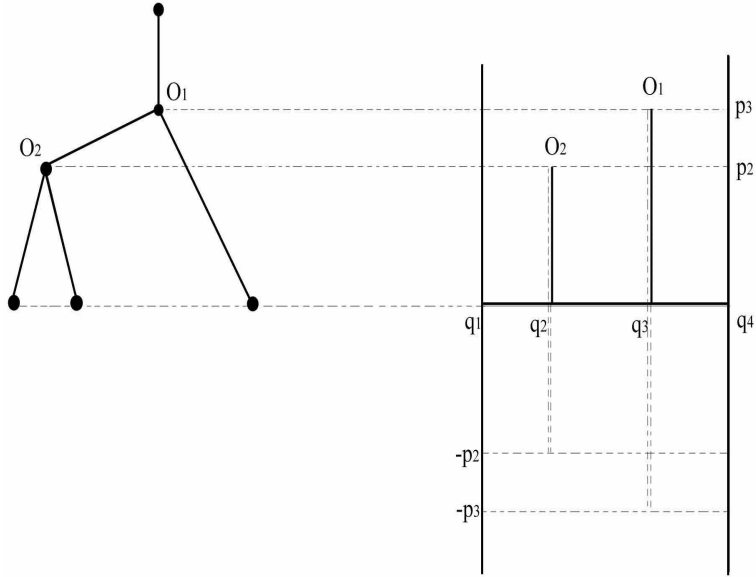}
\caption{1-dimensional model problem}
\end{figure}

Assume again now that the collisions with the walls are not
absolutely elastic: if the particle hits the wall at the point $q_k$
with a speed $p$, $|p|\leq p_k$, it is reflected from $q_k$ with the
speed $-p(1-\varepsilon c_k(p))$, where $c_k(\cdot)$ are defined as
before but now $(1-\ve c_k(p))$ represents the coefficient of
restitution for the wall $\{(q,p)\in \mathbb{R}^2: 0<p\leq p_k,
q=q_k\}$, and $0<\varepsilon << 1$. Then for each $\varepsilon
>0$, the velocity tends to zero as $t \rightarrow \infty$. On each velocity level $p_k$,
$p_k<p_0$, phase trajectory goes to one of the two "wells" separated
by $q_k$. Corresponding phase space $\sqcap=\{(q,p)\in \mathbb{R}^2,
q_1\leq q \leq q_n\}$ is shown in Fig.3 for $n=4$. Notice that
actually each interior wall in the phase space $\sqcap$ consists of
two sides (the broken line represents another side in Fig.3). For
simplicity we will identify the two sides in our investigation and
we assume that the coefficients $c_k(\cdot)$ is the same for both
faces of the wall $q=q_k$.

We will see that the "choice of the well" at each $q_k$ is very
sensitive to $\varepsilon$ as $\varepsilon \downarrow 0$, and the
behavior of the trajectory becomes, actually, stochastic in a
certain sense. Now we will give the exact meaning to this statement
and describe the limiting stochastic process.

Let us, first, consider the case when the model system has only two
trapping wells. An example of such a model is shown in Fig.4. A
particle with unit mass starts its motion from a point $x=(q_0,p_0)$
in the phase space $\sqcap$ and has instantaneous reflection each
time when it hits the boundary $q=-a_1$(corresponding to well
$\mathcal{E}_1$) or $q=a_2$(corresponding to well $\mathcal{E}_2$).
Each time when the particle hits the boundary $q=-a_1$ or $q=a_2$ it
will be reflected and move to the point $(-1,-p(1-\varepsilon
c_1(p)))$ or $(1,-p(1-\varepsilon c_2(p)))$, respectively. The
coefficient of restitution for the middle wall (notice that it has
two sides) $q=0$ is $1-\varepsilon c_3(\cdot)$. After the particle
enters one of the wells $\mathcal{E}_1$ or $\mathcal{E}_2$, it hits
corresponding boundaries and loses energy according to the
coefficient of restitution of the corresponding boundary. We denote
such a motion by
$X_t^{\ve,x}=X_t^\varepsilon=(q^\varepsilon(t),p^\varepsilon(t))$,
where $(q^\ve(0), p^\ve(0))=(q_0,p_0)=x$. The superscript $x$
represents the starting point of the motion. Here and henceforth,
when we use the notation $X^\ve_t$, it means that we have omitted
the superscript $x$. Put
$H^\varepsilon(t)=\displaystyle{\frac{(p^\varepsilon(t))^2}{2}}$. We
rescale time and let
$\widetilde{X}_t^\varepsilon=X_{t/\varepsilon}^\varepsilon$. As
before, we consider the slow component of $\widetilde{X}_t^\ve$
which is the projection
$Y_t^\varepsilon=(H(\widetilde{X}_t^\varepsilon),
K(\widetilde{X}_t^\varepsilon))$ of $\widetilde{X}_t^\ve$ onto the
graph $\Gamma$. Here
$H(\widetilde{X}_t^\varepsilon)=H^\varepsilon(t/\varepsilon)$. Since
the coefficient of restitution depends only on the magnitude of
velocity but not its direction, i.e. $c_i(-p)=c_i(p)$, $i=1,2$, we
can write $c_i(-p)=c_i(p)=c_i(\sqrt{2H})$. Now we prove:

\begin{figure}
\includegraphics[height=8cm, width=15cm , bb=73 32 563 303]{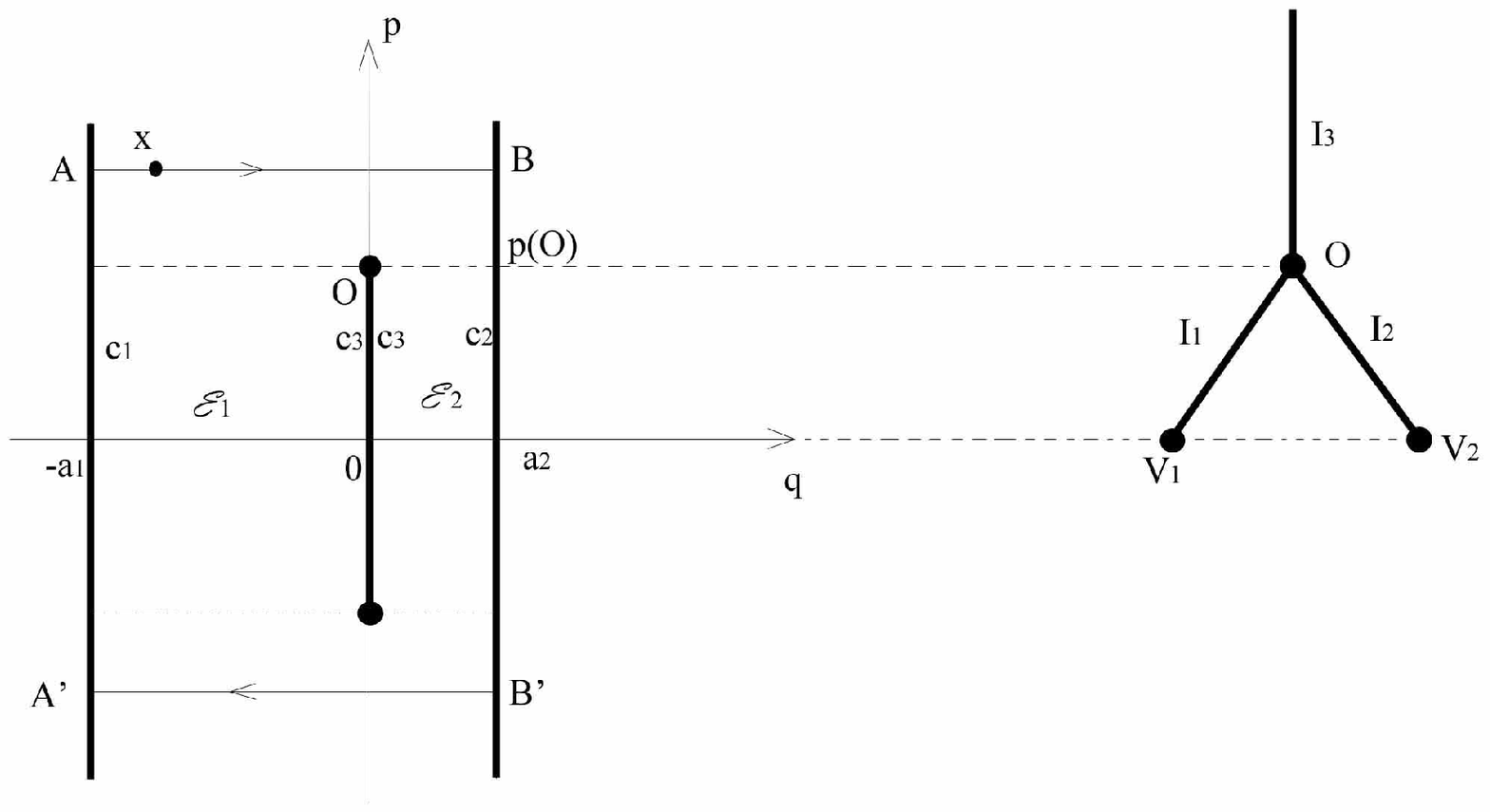}
\caption{The case when the model system has only two wells}
\end{figure}

\

\textbf{Lemma 2.1} \textit{Within each edge of the graph $\Gamma$,
as $\varepsilon \downarrow 0$, the first component of the process
$Y_t^\varepsilon$, i.e.
$H(\widetilde{X}_t^\varepsilon)=\displaystyle{\frac{(p^\varepsilon(t/\varepsilon))^2}{2}}$,
converges uniformly to a deterministic motion $H(t)$ which satisfies
the differential equation
$$\frac{dH}{dt}=-2\frac{c_1(\sqrt{2H})+c_2(\sqrt{2H})}{T_3(H)}H \ \ , \ \ H(0)=H_0 \
\ on \ \ I_3, \eqno(2.1)$$ and
$$\frac{dH}{dt}=-2\frac{c_1(\sqrt{2H})+c_3(\sqrt{2H})}{T_1(H)}H , \ \ on \ \
I_1,\eqno(2.2)$$
$$\frac{dH}{dt}=-2\frac{c_2(\sqrt{2H})+c_3(\sqrt{2H})}{T_2(H)}H , \
\ on \ \ I_2,\eqno(2.3)$$ respectively. Here
$T_3(H)=\displaystyle{\frac{2(a_1+a_2)}{\sqrt{2H}}}$,
$T_1(H)=\displaystyle{\frac{2a_1}{\sqrt{2H}}}$,
$T_2(H)=\play{\frac{2a_2}{\sqrt{2H}}}$ are the corresponding periods
of motion on phase picture $\sqcap$ along the energy level $H$ for
each edge of $\Gamma$.}

\

\textit{Proof}: Let us prove, for example, (2.1). The proofs of
(2.2) and (2.3) are exactly similar to the proof of (2.1). We prove
this result by using a slight modification of the standard method of
justification of the averaging principle (compare with [1], Ch.10).
First of all, we view the whole contour (see Fig.4) $A\ra B \ra B'
\ra A' \ra A$ on which the particle has a fixed amount of energy $H$
as a loop $G_H$ (i.e. identify $A$ with $A'$ and $B$ with $B'$). One
can introduce coordinates $(Q^\ve_{t/\ve}, H^\ve_{t/\ve})$ to
describe the motion of the particle on this circle. Here
$Q^\ve_{t/\ve}=q^\ve_{t/\ve}+2(k-1)(a_1+a_2)$ if the particle is the
$k$-th time on the upper half plane of the phase space, and
$Q^\ve_{t/\ve}=-q^\ve_{t/\ve}+2(k-1)a_1+2ka_2$ if the particle is
the $k$-th time on the lower half plane of the phase space. We
identify $Q^\ve_{t/\ve}$ with $Q^\ve_{t/\ve}+2k(a_1+a_2)$, $k \in
\mathbb{Z}$. The motion $(Q^\ve_{t/\ve}, H^\ve_{t/\ve})$ has
 a fast component $Q^\ve_{t/\ve}$ and a slow component $H^\ve_{t/\ve}$
and they satisfy the following equation

\begin{equation*}
\left\{\begin{array}{lr}
\dot{Q}^\ve_{t/\ve}=\play{\frac{1}{\ve}\sqrt{2H^\ve_{t/\ve}}} \ ,\\
\dot{H}^\ve_{t/\ve}=-[\Delta(Q^\ve_{t/\ve}-a_2)+\Delta(Q^\ve_{t/\ve}+a_1)](2c(Q^\ve_{t/\ve},\sqrt{2H^\ve_{t/\ve}})-\ve
c^2(Q^\ve_{t/\ve}, \sqrt{2H^\ve_{t/\ve}}))H^\ve_{t/\ve} \ .
\end{array}\right.\eqno(2.4)
\end{equation*}

Here $\Delta(\cdot)$ is the Dirac $\delta$-function and
$c(Q^\ve_{t/\ve}, \sqrt{2H^\ve_{t/\ve}})=c_1(\sqrt{2H^\ve_{t/\ve}})$
if $Q^\ve_{t/\ve}=-a_1$ and $c(Q^\ve_{t/\ve},
\sqrt{2H^\ve_{t/\ve}})=c_2(\sqrt{2H^\ve_{t/\ve}})$ if
$Q^\ve_{t/\ve}=a_2$.

Now we introduce an auxiliary function $u=u(Q,H)$ satisfying the
following differential equation

$$\sqrt{2H}\frac{\partial u}{\partial Q}= A-\overline{A} \ , \ Q\in G_H \eqno(2.5)$$

where $$A=-[\Delta(Q-a_2)+\Delta(Q+a_1)](2c(Q,\sqrt{2H})-\ve c^2(Q,
\sqrt{2H}))H,$$

and
$$\overline{A}=-2\play{\frac{c_1(\sqrt{2H})+c_2(\sqrt{2H})-\ve(c_1^2(\sqrt{2H})+c_2^2(\sqrt{2H}))}{T_3(H)}}$$
is the average of $A$ along the loop $A\ra B\ra B' \ra A' \ra A$.

One can check the following properties of $u(Q,H)$: (i) the
differential equation (2.5) describing $u$ really has a well-defined
solution $u(Q,H)$ on the loop $G_H$; (ii) the solution $u(Q,H)$ is
uniformly bounded together with its first derivatives for $Q\in
G_H$, $H \in [H(O), M]$ for any $M\in (H(O), \infty)$.

Now, by the fundamental theorem of calculus, we have:

$$\begin{array}{l}
u(Q^\ve_{t/\ve}, H^\ve_{t/\ve})-u(Q^\ve_0, H^\ve_0)\\
\play{=\int_0^t\frac{\partial u}{\partial
Q}\dot{Q}^\ve_{s/\ve}ds+\int_0^t\frac{\partial u}{\partial
H}\dot{H}^\ve_{s/\ve}ds}\\
\play{=\frac{1}{\ve}\int_0^t(A-\overline{A})ds+O(1)}.
\end{array}$$

Therefore we see that $\play{\max\limits_{0\leq t \leq
T}\left|\int_0^t(A-\overline{A})ds\right|\leq C \ve}$ for some
positive $C>0$.

Now, we compare the two differential equations:

$$\dot{H}^\ve_{t/\ve}=-[\Delta(Q^\ve_{t/\ve}-a_2)+\Delta(Q^\ve_{t/\ve}+a_1)](2c(Q^\ve_{t/\ve},\sqrt{2H^\ve_{t/\ve}})-\ve
c^2(Q^\ve_{t/\ve}, \sqrt{2H^\ve_{t/\ve}}))H^\ve_{t/\ve}\ , \ \
H^\ve_0=H_0;$$

$$\frac{dH}{dt}=-2\frac{c_1(\sqrt{2H})+c_2(\sqrt{2H})}{T_3(H)}H \ , \ \ H(0)=H_0.$$

Let $m(T)=\max\limits_{0\leq t \leq T}|H^\ve_{t/\ve}-H(t)|$. We
have, that

$$\begin{array} {l} |H^\ve_{t/\ve}-H(t)|\\
=\play{\left|\int_0^t (AH^\ve_{s/\ve}-\overline{A}H(s))ds\right| +
2\ve B}\\
\leq\play{\int_0^t
|A-\overline{A}|H^\ve_{s/\ve}ds+\int_0^t\overline{A}|H^\ve_{s/\ve}-H(s)|ds+2\ve
B}.
\end{array}$$

The last inequality implies that $\play{m(T)\leq (CD+2B)\ve+ E
\int_0^Tm(s)ds}$, where $B$ is the uniform bound for
$\play{\frac{c_1^2(\sqrt{2H})+c_2^2(\sqrt{2H})}{T_3(H)}}$, $D$ is
the uniform bound for $H^\ve_{s/\ve}$, and $E$ is the uniform bound
for $\overline{A}$. By using Gronwall's inequality, we see that
$m(T)\leq (CD+2B)\ve e^{ET}$, i.e. $\lim\limits_{\ve \da
0}\max\limits_{0\leq t \leq T}|H^\ve_{t/\ve}-H(t)|=0$.
$\blacksquare$

\

\

Notice that in this case $T_i(H)\sim H^{-\frac{1}{2}}$ as $H
\rightarrow H(O)$, this means that $Y^\varepsilon_t$ will enter the
interior vertex $O$ of the graph $\Gamma$ in finite time. Thus one
might ask what is the behavior of the motion $Y^\varepsilon_t$ at
the interior vertex $O$. First, we regularize this problem by
considering small stochastic perturbation of the initial conditions
in our model.

\

\textbf{Lemma 2.2} \textit{Let $\mathcal{U}(x,\delta)=\{y\in \sqcap:
|x-y|<\delta\}$, $x\in \sqcap$, $\delta > 0$,
$\mathcal{U}_i^\ve=\{x\in \sqcap: X_\cdot^{\ve,x}$ eventually enters
the well $\mathcal{E}_i\}$, $i\in \{1,2\}$ (see Fig.5), where
$X_\cdot^{\ve,x}$ is the perturbed trajectory starting at
$X_0^{\ve,x}=x\in \sqcap$. Assume that $H_0=H(x)>H(O)$. Then}

\begin{equation*}
\lim\limits_{\delta \downarrow 0}\lim\limits_{\varepsilon \downarrow
0}\frac{\mu(\mathcal{U}^\varepsilon_1\bigcap\mathcal{U}(x,\delta))}{\mu(\mathcal{U}^\varepsilon_2\bigcap\mathcal{U}(x,\delta))}=\frac{c_1(H(O))}{c_2(H(O))},
\end{equation*}

\textit{where $\mu$ is the Lebesgue measure in $\mathbb{R}^2$.}

\

\textit{Proof}: Without loss of generality one can assume that
$p_0>0$. We see, from Fig.5, that $\mathcal{U}(x,\delta)$ is covered
by narrow shaded or white strips, where shaded strips belong to
$\mathcal{U}^\varepsilon_1$, and white strips belong to
$\mathcal{U}^\varepsilon_2$. For a fixed $x$, consider the nearest
to $x$ shaded and white strip (ordered as a pair, shaded on top).
The upper $p$-level of the shaded strip of the pair is denoted by
$a$. And the lower $p$-level of the shaded strip (also the upper
$p$-level of the white strip in the pair) by $b$. Let the lower
$p$-level of the white strip in the pair be $c$. Let $n=n(\ve,x)$ be
the number of collisions that the particle made with either of the
walls $q=-a_1$ or $q=a_2$ before entering one of the wells.

Let us first estimate the $p$-width of the shaded and white strips.
Put
$$A(a)=a(1-\varepsilon c_1(a))(1-\varepsilon c_2(a(1-\varepsilon
c_1(a))))=a(1-\varepsilon f_1^\varepsilon(a)+\varepsilon^2
f_2^\varepsilon(a))$$ where
$f_1^\varepsilon(a)=c_1(a)+c_2(1-\varepsilon c_1(a))$ and
$f_2^\varepsilon(a)=c_1(a)c_2(1-\varepsilon c_1(a))$. Then the
boundaries of the next shaded and white strip has $p$-height $A(a)=c
, A(b), A(c)$, respectively.

We see that
$$|c-a|=|A(a)-a|=\varepsilon a|f_1^\varepsilon (a)-\varepsilon f_2^\varepsilon (a)|,$$

So that
$$M_1\varepsilon\leq|c-a|\leq
M_2\varepsilon$$ for some $M_1, M_2>0$. Therefore,
$n=n(x,\varepsilon)\sim O(\displaystyle{\frac{1}{\varepsilon}})$.

\begin{figure}
\centering
\includegraphics[height=7cm, width=9cm , bb=143 10 429 306]{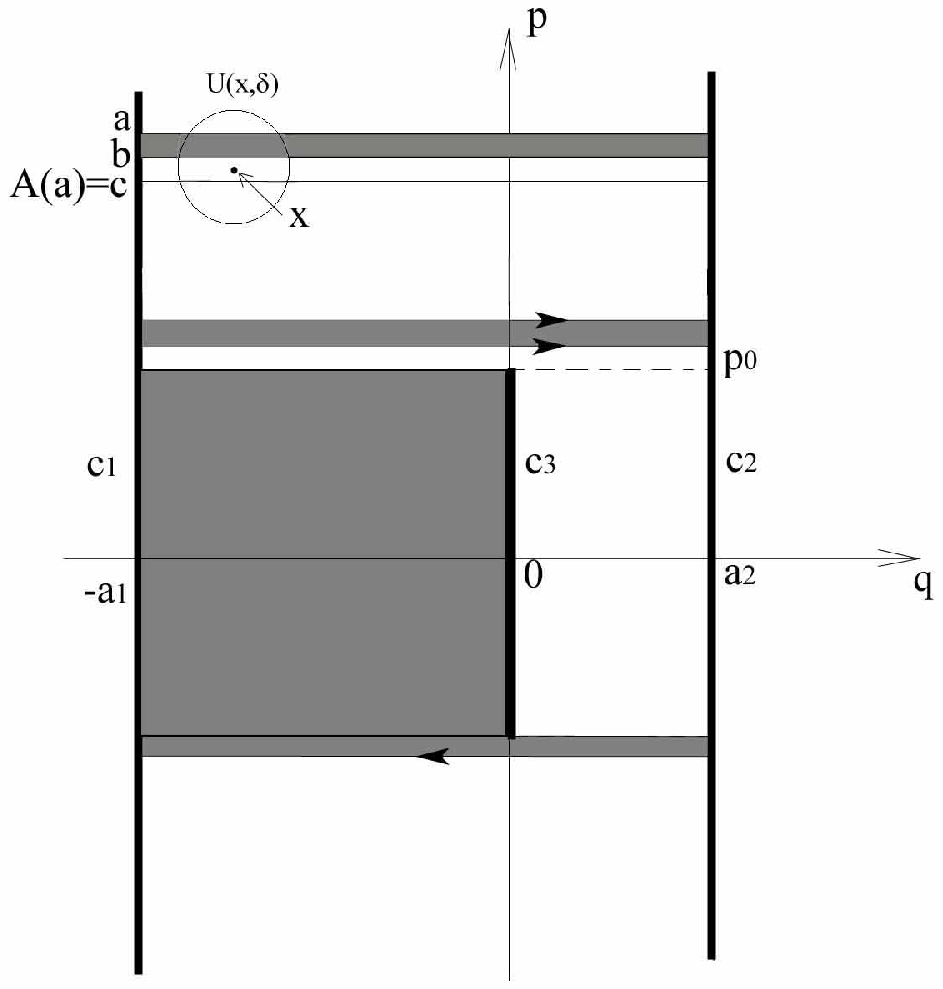}
\caption{Regularization by stochastic perturbation of the initial
condition}
\end{figure}

We have

$$\begin{array}{l}
A(a)-A(b)
\\=a-b-\varepsilon(af_1^\varepsilon(a)-bf_1^\varepsilon(b))+\varepsilon^2(af_2^\varepsilon(a)-bf_2^\varepsilon(b))
\\=a-b-\varepsilon(u_1^\varepsilon(b)(a-b)+u_2^\varepsilon(\xi)(a-b)^2)+\varepsilon^2(v_1^\varepsilon(b)(a-b)+v_2^\varepsilon(\eta)(a-b)^2)
\\=(a-b)(1-\varepsilon u_1^\varepsilon(b)+O(\varepsilon^2))
\end{array}
$$

where $u_i$, $v_i$, $i\in \{1,2\}$ are bounded smooth functions, and
$\xi, \eta \in [b,a]$.

 Similarly, $$A(b)-A(c)=(b-c)(1-\varepsilon
u_1^\varepsilon(c)+O(\varepsilon^2)).$$

These inequalities imply that

$$\begin{array}{l}
\displaystyle{\frac{A(a)-A(b)}{A(b)-A(c)}}
\\=\displaystyle{\frac{a-b}{b-c}\left(\frac{1-\varepsilon u_1^\varepsilon(b)+O(\varepsilon^2)}{1-\varepsilon
u_1^\varepsilon(c)+O(\varepsilon^2)}\right)}
\\=\displaystyle{\frac{a-b}{b-c}\left(1-\frac{\varepsilon(u_1^\varepsilon(b)-u_1^\varepsilon(c))+O(\varepsilon^2)}{1-\varepsilon
u_1^\varepsilon(c)+O(\varepsilon^2)}\right)}
\\=\displaystyle{\frac{a-b}{b-c}(1+O(\varepsilon^2))}.
\end{array}$$

Thus after $n=n(\varepsilon,x)\sim
\displaystyle{O(\frac{1}{\varepsilon})}$ steps the change of ratio
of the $p$-width of the nearby shaded and white strips is
asymptotically $(1+O(\varepsilon^2))^{n(\varepsilon,x)}\sim 1$, and
our lemma then immediately follows from the definition of the sets
$\mathcal{U}(x,\delta)$, $\mathcal{U}^\varepsilon_1$ and
$\mathcal{U}^\varepsilon_2$. $\blacksquare$

\

\

In order to include our result in the framework of weak convergence
of stochastic processes in the space of continuous functions, we
still need the following technical construction. Notice that, for
fixed $\varepsilon > 0$, the function $H^\varepsilon(t/\varepsilon)$
as a function of $t$, is a step function. We can assume that it is
continuous from the right and having limit from the left. We now
connect neighboring discontinuity points of the graph of
$H^\varepsilon(t/\varepsilon)$ by straight line segments. The
resulting trajectory is denoted by $\widehat{H}^\varepsilon(t)$. For
fixed $T>0$, we see that $\widehat{H}^\varepsilon(t)\in C_{0T}$. For
a positive constant $C>0$ and all $0<t<T$, we have
$$|H^\varepsilon(t/\varepsilon)-\widehat{H}^\varepsilon(t)|<C\varepsilon. \eqno(2.6)$$

Since the time distance between two discontinuity points is bounded
from below by the inverse of the initial velocity, the slope of all
the line segments of the graph of $\widehat{H}^\varepsilon(t)$ is
bounded uniformly in $\ve \in (0,1]$. Thus the family of functions
$\{\widehat{H}^\varepsilon(t)\}_{\varepsilon
>0}$ is an equicontinuous family in $C_{0T}$. Also, it is uniformly
bounded for $0<t<T$ and $\varepsilon\in (0,1]$. Thus by
Ascoli-Arzela Theorem the family
$\{\widehat{H}^\varepsilon(t)\}_{\varepsilon
>0}$ is compact in $C_{0T}$ with uniform topology.

\

Let $\xi_\delta$ be a two-dimensional random variable with a
continuous density $f_\delta(x)$ such that $f_\delta(x)>0$ for
$|x|<\delta$ and $f_\delta(x)=0$ for $|x|\geq \delta$. Let
$X_t^{\ve, \delta}=(p_t^{\ve, \delta}, q_t^{\ve, \delta})$ be the
trajectory of the nearly-elastic motion starting at $X_0^{\ve
,\delta}=x+\xi_\delta$ (the point $x=(q_0,p_0)\in \sqcap$ is fixed,
and we excluded it from the notation). Put $\widetilde{X}_t^{\ve,
\delta}=X^{\ve,\delta}_{t/\ve}$.

Consider a stochastic process $Y_t=(H(t),K(t))$, $0<t<T$, on the
graph $\Gamma$ defined by the conditions: $Y_0=Y(q_0,p_0)$; $H(t)$
is a deterministic motion inside each edge of $\Gamma$ satisfying
equations (2.1)-(2.3) respectively; if $Y_0=Y(q_0,p_0)\in I_3$, the
trajectory $Y_t$ reaches the interior vertex $O\in \Gamma$ in a
finite time, instantaneously leaves $O$ and enters the edges $I_1$
or $I_2$ with probabilities
$\play{p_1=\frac{c_1(H(O))}{c_1(H(O))+c_2(H(O))}}$,
$\play{p_1=\frac{c_2(H(O))}{c_1(H(O))+c_2(H(O))}}$, respectively.
The conditions listed above define the process $Y_t$ in a unique way
(in the sense of distributions).

Taking into account (2.6), the remark concerning compactness and
Lemmas 2.1 and 2.2, we get the following

\

\textbf{Theorem 2.1} \textit{The process
$\widehat{Y}_t^{\ve,\delta}=(\widehat{H}(\widetilde{X}^{\ve,\delta}_t),K(\widetilde{X}^{\ve,\delta}_t))$
on $\Gamma$ converges weakly in the space of continuous functions
$[0,T]\ra \Gamma$ equipped with the uniform topology to the process
$Y_t=(H(t),K(t))$ when, first $\ve \da 0$ and then $\delta \da 0$.}

\

\

We now discuss the case when the model problem has more than two
energy trapping wells. It turns out that in the multi-well case, the
perturbation of the initial condition, in general, will \textit{not}
lead to a regularization of the problem: the limit of
 $Y_t^{\ve,\delta}$ as $\ve \da 0$ may not exist.

Consider an example as shown in Fig.6. There are 3 trapping wells,
1, 2 and 3. Wells 2 and 3 are separated on a $p$-level $p^*_b$ and
they combined together are separated from well 1 on a $p$-level
$p^*_a$. $p^*_a>p^*_b$. We call the wells 2 and 3 combined together
well 4. Well 1 and well 4 combined together is called well 5.
Suppose the final strip $A$ as drawn in the picture is for those
initial values that finally enter well 4. The boundary restitution
coefficients are positive constants $c_1$ and
$c_2=c_3=c_4=c_5=c(\neq c_1)$ as shown in Fig.6.

For each $\ve > 0$, the shadowed strip $A$ on the top of wells 1 and
4 is moved by the nearly-elastic dynamics to the well 4. The strip
$A$ has the width $c \ve$. The well 4 between $p$-levels $p_b^*$ and
$p_a^*$ is also covered by strips of width $c \ve$ (suppose that
$\ve$ is such that an integer number of such strips is situated
between $p_b^*$ and $p_a^*$). The highest strip $B$ in the wall 4 is
moved by our dynamics (as a whole) alternatively either to well 2 or
to well 3 as $\ve \da 0$. But the strip $A$ in one step (one
reflection from the right exterior wall) is going to the strip $B$.
Therefore, for the sequence of $\ve \da 0$ for which
$\play{\frac{p_a^*-p_b^*}{\ve}}$ is an integer, the nearly-elastic
trajectory starting from any $x=(q_0,p_0)$, $p_0>p_a^*$, for $t$
large enough is distributed alternatively either between walls 1 and
2 or between walls 1 and 3. This means that the limiting
distribution for $X_t^{\ve, \delta}$ as $\ve \da 0$ does not exist,
and the problem cannot be regularized by random perturbations of the
initial point.

\begin{figure}
\centering
\includegraphics[height=8cm, width=15cm , bb=44 73 514 239]{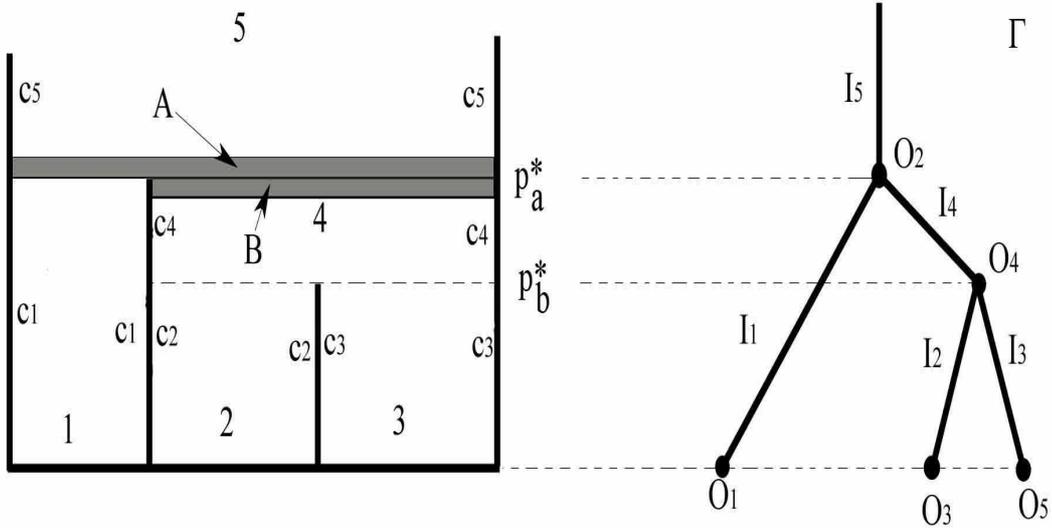}
\caption{The case when the system has more than two trapping wells}
\end{figure}

\section{Regularization by stochastic perturbation of the dynamics}

Let us now consider small stochastic perturbation of the dynamics
rather than the initial condition. We will show that, this
regularization works when the system has any number of trapping
wells. But let us first start from the case of system with only two
trapping wells.

The particle of unit mass starts its motion from $x=(q_0,p_0)$ with
$p_0>0$ (see Fig.7). Each time when the particle hits the wall
$q=-a_1$ or $q=a_2$, it instantaneously moves to
$(-a_1,-p(1-\varepsilon c_1-\varepsilon \delta \eta_k))$, $p<0$ and
$(a_2,-p(1-\varepsilon c_2-\varepsilon \delta \xi_k))$, $p>0$,
respectively, and then it is reflected. We assume for simplicity
that $c_1, c_2$ are positive constants, $0<\varepsilon<<1$ and
$0<\delta<<1$. The restitution coefficient for the interior wall is
assumed to be $1-\varepsilon c_3-\varepsilon \delta \zeta_k$ for
both faces. And $c_3$ is also a positive constant.

Here $\xi_k$, $\eta_k$, $\zeta_k$, $k=1,2,...$, are independent
sequences of i.i.d random variables with continuous densities; the
random variable with subscript $k$ is used in the restitution
coefficient when the trajectory hits the corresponding wall $k$-th
time. In the following we also use the general notation $\xi$,
$\eta$, $\zeta$ to denote independent r.v.'s which has the same
distribution as $\xi_k$, $\eta_k$ and $\zeta_k$, respectively. We
assume for brevity that random variables $\xi$, $\eta$ and $\zeta$
are bounded with probability 1. Then, without loss of generality,
one can consider just strictly positive $\xi$, $\eta$ and $\zeta$:
$P\{\alpha < \xi < \beta\}=P\{\alpha < \eta < \beta\}=P\{\alpha <
\zeta < \beta\}=1$ for some $0<\alpha < \beta < \infty$. Actually,
one can replace the boundedness of random variables by their
positivity and the finiteness of their moments (together with the
existence of a continuous density).

The resulting process is denoted by
$X^{\varepsilon,\delta}_t=(q^{\varepsilon,\delta}(t),p^{\varepsilon,\delta}(t))$.
We rescale time and put
$\widetilde{X}^{\varepsilon,\delta}_t=X^{\varepsilon,\delta}_{t/\varepsilon}$.
Let
$Y_t^{\varepsilon,\delta}=(H(\widetilde{X}_t^{\varepsilon,\delta}),K(\widetilde{X}_t^{\varepsilon,\delta}))$
be the projection of $\widetilde{X}_t^{\ve, \delta}$ onto the graph
$\Gamma$ (see Fig.7). Put
$H^{\varepsilon,\delta}(t)=H(\widetilde{X}_t^{\varepsilon,\delta})$.
For fixed $\varepsilon >0$, this is a stochastic process with jumps
and each sample path is a step function which is right continuous
and has limit from the left. Using the same construction as in the
end of Section 2, we obtain from $H^{\ve, \delta}(t)$ a continuous
piecewise linear process $\widehat{H}^{\varepsilon, \delta}(t)$. We
also put
$\widehat{Y}_t^{\varepsilon,\delta}=(\widehat{H}(\widetilde{X}_t^{\varepsilon,\delta}),K(\widetilde{X}_t^{\varepsilon,\delta}))$.
Since $\xi$ and $\eta$ are bounded,
$$P(|\widehat{H}^{\varepsilon, \delta}(t)-H^{\varepsilon,\delta}(t)|<C\varepsilon)=1$$
for some constant $C>0$.

Now we are going to prove the weak convergence of
$\widehat{Y}_t^{\varepsilon,\delta}$ in the space $C_{0T}(\Gamma)$
of continuous functions $[0,T]\ra \Gamma$ provided with uniform
topology, as first $\varepsilon \downarrow 0$ and then $\delta
\downarrow 0$, to a stochastic process $Y_t$ on the graph $\Gamma$.
We do this through a series of lemmas.

\begin{figure}
\includegraphics[height=8cm, width=15cm , bb=65 30 559 304]{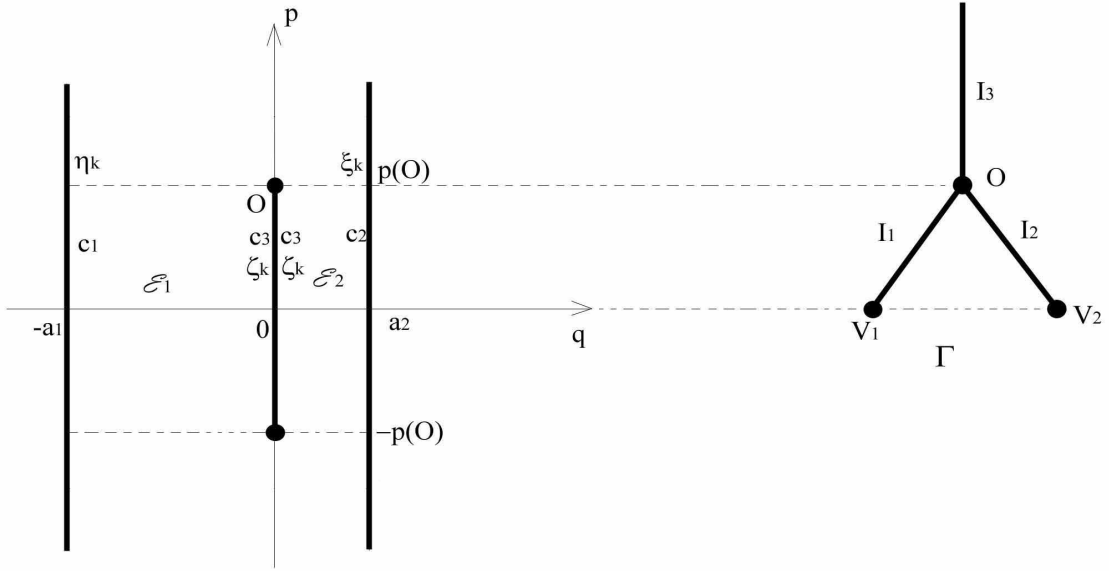}
\caption{Regularization by stochastic perturbation of the dynamics}
\end{figure}

\

\textbf{Lemma 3.1} \textit{For fixed $\delta >0$, the family
$\{\widehat{Y}^{\varepsilon, \delta}(t)\}_{\varepsilon >0}$ is tight
in the space $C_{0T}(\Gamma)$}.

\

\textit{Proof}: We have to verify the following (see [7], Ch.6,
Theorem 6.4.2): there exist $\alpha, \beta > 0,$ s.t. for any $h >
0$, any $0< t < t+h <T$, any $\varepsilon > 0$,

$$\mathbb{E}\left(\rho(\widehat{Y}_t^\varepsilon,\widehat{Y}_{t+h}^\varepsilon)\right)^\alpha \leq Mh^{1+\beta}.$$

Here $\rho(x,y)$ denotes the distance between the points $x,y\in
\Gamma$ defined as follows: if $x$ and $y$ lies in the same edge
$I_i$ ($i=1,2,3$) of the graph $\Gamma$, then $\rho(x,y)$ is just
the usual Euclidean distance between points $x$ and $y$ within the
line segment; otherwise, suppose $x\in I_i$ and $y\in I_j$, $i\neq
j$, then consider a path lying on the graph $\Gamma$ connecting $x$
and $y$, we denote $\rho(x,y)$ be just the minimal length of such
paths. Since our $c_1$, $c_2$, $c_3$, and $\xi, \eta, \zeta$ are
bounded with probability 1, and the period $T(H)$ of the
non-perturbed (elastic) motion is bounded for $0<H<\infty$ and
separated from zero, we see that
$P(\rho(\widehat{Y}_t^\varepsilon,\widehat{Y}_{t+h}^\varepsilon)<Mh)=1$
 for some constant $M>0$. Therefore one can take $\alpha=2, \beta=1$, and the statement
follows.$\blacksquare$

\

\textbf{Lemma 3.2} \textit{Let $H^{\ve,\delta}(0)=H_0>H(O)$. Within
each edge of the graph $\Gamma$, as $\varepsilon \downarrow 0$, the
process $\widehat{H}^{\varepsilon,\delta}(t)$, for $0<t<T<\infty$,
converges uniformly in probability to a deterministic motion
$H^\delta(t)$ which is defined by the equations
$$H^\delta(t)=\left(\frac{t}{\sqrt{2}}\frac{c_1+c_2+\delta (\mathbb{E}\xi+\mathbb{E}\eta)}{a_1+a_2}+H_0^{-1/2}\right)^{-2}, \ on \ I_3 ,$$ and
$$H^\delta(t)=\left(\frac{t-t_0}{\sqrt{2}}\frac{c_1+c_3+\delta (\mathbb{E}\zeta+\mathbb{E}\eta)}{a_1}+H(O)^{-1/2}\right)^{-2}, \ on \ I_1 ,$$
$$H^\delta(t)=\left(\frac{t-t_0}{\sqrt{2}}\frac{c_2+c_3+\delta (\mathbb{E}\xi+\mathbb{E}\zeta)}{a_2}+H(O)^{-1/2}\right)^{-2}, \ on \ I_2, $$
respectively. Here $H(O)$ is the energy corresponding to the
interior vertex $O$ and
$t_0^\delta=\play{\frac{\sqrt{2}(a_1+a_2)(H(O)^{-1/2}-H_0^{-1/2})}{c_1+c_2+\delta
\mathbb{E}(\xi+\eta)}}$ is the time for the motion $H^\delta(t)$ to
come to the interior vertex $O$.}

\

\textit{Proof}: The proof of this lemma is similar to the proof of
Lemma 2.1 and we use the same notations. The system (2.4) should be
replaced by the following system:

$$
\left\{\begin{array}{lr}
\dot{Q}^\ve_{t/\ve}=\play{\frac{1}{\ve}\sqrt{2H^\ve_{t/\ve}}} \ ,\\
\dot{H}^\ve_{t/\ve}=-[\Delta(Q^\ve_{t/\ve}-a_2)+\Delta(Q^\ve_{t/\ve}+a_1)](2c(Q^\ve_{t/\ve})+2\delta\cdot
r(Q^\ve_{t/\ve})-\ve b(Q^\ve_{t/\ve}))H^\ve_{t/\ve} \ .
\end{array}\right.$$ $$\eqno(3.1)
$$

Here $\Delta(\cdot)$ is the Dirac $\delta$-function so that the
right hand side of the last equation in (3.1) is not zero just if
$Q^\ve_{t/\ve}=a_1$ or $Q^\ve_{t/\ve}=a_2$; $r(Q^\ve_{t/\ve})=\xi_k$
and $c(Q^\ve_{t/\ve})=c_2$ when $c(Q^\ve_{t/\ve})=a_2$;
$r(Q^\ve_{t/\ve})=\eta_k$ and $c(Q^\ve_{t/\ve})=c_1$ when
$c(Q^\ve_{t/\ve})=a_1$. The function
$b(Q^\ve_{t/\ve})=\left(c(Q^\ve_{t/\ve})+\delta\cdot
r(Q^\ve_{t/\ve})\right)^2$ is uniformly bounded for $0\leq t \leq T$
with probability one. After that, using the same arguments as in
Lemma 2.1 and the law of large numbers, we obtain an equation for
the limiting slow component $H(t)$ on each edge. For constant $c_i$,
these equations can be solved explicitly and we get the statement of
Lemma 3.2. $\blacksquare$

\

Lemma 3.2 implies the following:

\

\textbf{Corollary 3.1} \textit{Within each edge of the graph
$\Gamma$, for the time $0<t<T<\infty$, the process
$\widehat{H}^{\varepsilon,\delta}(t)$, as $\varepsilon , \delta
\downarrow 0$, converges in probability to a deterministic motion
$H(t)$ defined by the equations
$$H(t)=\left(\frac{t}{\sqrt{2}}\frac{c_1+c_2}{a_1+a_2}+H_0^{-1/2}\right)^{-2}, \ on \ I_3 \eqno(3.2)$$ and
$$H(t)=\left(\frac{t-t_0}{\sqrt{2}}\frac{c_1+c_3}{a_1}+H(O)^{-1/2}\right)^{-2}, \ on \ I_1 \eqno(3.3)$$
$$H(t)=\left(\frac{t-t_0}{\sqrt{2}}\frac{c_2+c_3}{a_2}+H(O)^{-1/2}\right)^{-2}, \ on \ I_2 \eqno(3.4)$$
respectively. Here
$t_0=\play{\frac{\sqrt{2}(a_1+a_2)(H(O)^{-1/2}-H_0^{-1/2})}{c_1+c_2}}$.}

\

Let us consider now the slow motion near the interior vertex $O$. We
first prove the following auxiliary lemma concerning random walks.
Let $\{\xi_k\}$ and $\{\eta_k\}$ be independent sequences of i.i.d
random variables. Assume that the random variables have continuous
densities and $P\{\alpha<\xi<\beta\}=P\{\alpha<\eta<\beta\}=1$ for
some $0<\alpha<\beta<\infty$. Put
$$S_0^x=x \ , \ S_{2m}^x=x+\sum\li_{k=1}^{m}(\xi_k+\eta_k) \ , \ S_{2m+1}^x=x+S_{2m}^x+\xi_{m+1} \ .$$
It is clear that $S_n^x=x+S_n^0$, $P\{n\alpha<S_n^0<n\beta\}=1$.

Define $\tau_n^{x,\lambda}$ as the first time $m$ when $S_m^x$ is
greater than $n\lambda$: $\tau_n^{x,\lambda}=\min\{m:
S_m^x>n\lambda\}$.

Since $\mathbb{E}(\xi_k+\eta_k)>0$, the law of large numbers implies
that $P\{\tau_n^{x,\lambda}<\infty\}=1$ for any $x\in \mathbb{R}^1$,
$\lambda>0$, $n\in \mathbb{Z}$.

We use two equivalent types of notations: The initial point $x$ of
the random walk will be included either as a superscript, like
$S_n^{x}$, $\tau_n^{x,\lambda}$, or as a subscript in probabilities
and expected values so that $\mathbb{E}_xf(S_n)\equiv
\mathbb{E}f(S_n^x)$, $P_x\{\tau_n^\lambda<t\}\equiv
P\{\tau_n^{x,\lambda}<t\}$.

\

\

\textbf{Lemma 3.3} \textit{Under the conditions mentioned above,
$\lim\li_{n \ra \infty}P_0\{\tau_n^\lambda \ \ is \ \
even\}=\play{\frac{\mathbb{E}\eta}{\mathbb{E}\xi+\mathbb{E}\eta}}$,
$\lim\li_{n \ra \infty}P_0\{\tau_n^\lambda \ \ is \ \
odd\}=\play{\frac{\mathbb{E}\xi}{\mathbb{E}\xi+\mathbb{E}\eta}}$. }
\

\

\

\textit{Proof:} Put $m_1(n)=m_1=\play{\left[\frac{n
\lambda}{2\beta}\right]}$ if the latter integer is even, and
$m_1(n)=m_1=\play{\left[\frac{n \lambda}{2\beta}\right]-1}$
otherwise. It is clear that $S_{m_1}^0<\play{\frac{n\lambda}{2}}$.
Put $M=\mathbb{E}(\xi_k+\eta_k)$, $D=Var(\xi_k+\eta_k)$. Let
$\widetilde{m_1}(n)=\widetilde{m_1}=\dfrac{m_1}{2}$. Let
$f_{m_1}(x)$ be the density of $S_{m_1}^0-\widetilde{m_1}M$.

It follows from the central limit theorem that for each $\delta>0$
one can choose $N=N(\delta)$ so that

$$P\{-N\sqrt{\widetilde{m_1}(n)}<S^0_{m_1}-\widetilde{m_1}M<N\sqrt{\widetilde{m_1}(n)}\}>1-\delta \eqno(3.5)$$
for all large enough $n$.

Using the Markov property of the random walk $S_n^x$, we get

$$P_0\{\tau_n^\lambda \ is \ even \}=\mathbb{E}_0P_{S_{m_1(n)}}\{\tau_n^\lambda \ is \ even \}. \eqno(3.6)$$

One can conclude from (3.5) and (3.6) that

$$\left|P_0\{\tau_n^\lambda \ is \ even \}-
\int_{-N\sqrt{\widetilde{m_1}}}^{N\sqrt{\widetilde{m_1}}}
f_{m_1}(x)P_{x+\widetilde{m_1}M}\{\tau_n^\lambda \ is \ even
\}dx\right|<\delta. \eqno(3.7)$$

Since our random variables are bounded and have density, one can
apply the local central limit theorem to $S_{m_1(n)}^0$:

$$\left|f_{m_1}(x)-
\frac{1}{\sqrt{2\pi \widetilde{m_1}
D}}\exp(-\frac{x^2}{2\widetilde{m_1}
D})\right|<\frac{\rho^{(1)}_{m_1}}{\sqrt{\widetilde{m_1}}}
\eqno(3.8)$$ uniformly in $x\in \mathbb{R}^1$; here and later, we
denote by $\rho_{\ell}^{(k)}$ such sequences that $\lim\li_{\ell\ra
\infty}\rho_{\ell}^{(k)}=0$.

Divide the interval $\left[-N\sqrt{\widetilde{m_1}(n)},
N\sqrt{\widetilde{m_1}(n)}\right]$ into
$\left[\sqrt[4]{\widetilde{m_1}(n)}\right]$ equal intervals $I_1$,
..., $I_r$. It is clear that $r=r(n)\sim n^{1/4}$ and the length
$|I_k|$ of each interval $I_k$ is of order $n^{1/4}$ as $n \ra
\infty$. Because of (3.7) and (3.8),

$$\left|P_0\{\tau_n^\lambda \ is \ even \}-\sum\li_{k=1}^{r(n)}\int_{I_k}
\exp(-\frac{x^2}{2\widetilde{m_1}D})\frac{P_{x+\widetilde{m_1}M}\{\tau_n^\lambda
\ is \ even \} dx}{\sqrt{2\pi \widetilde{m_1} D}}\right|<\delta +
\rho_n^{(2)}. \eqno(3.9)$$

On each $I_k$, $k\in \{1,...,r\}$, choose a point $z_k$ such that
$$\int_{I_k}\exp(-\frac{x^2}{2\widetilde{m_1} D})dx=|I_k|\exp(-\frac{z_k^2}{2\widetilde{m_1} D}).$$

If $x\in I_k$, then
$$\left|\exp(-\frac{x^2}{2\widetilde{m_1}D})-\exp(-\frac{z_k^2}{2\widetilde{m_1}D})\right|
\leq \rho_n^{(3)}\exp(-\frac{z_k^2}{2\widetilde{m_1}D}),
\eqno(3.10)$$ where $\lim\li_{n \ra \infty}\rho_n^{(3)}=0$ since
$|x|<N\sqrt{\widetilde{m_1}(n)}$,
$|z_k|<N\sqrt{\widetilde{m_1}(n)}$, $|x-z_k|<\text{Const}\cdot n
^{1/4}$, $\widetilde{m_1}(n)$ is of order $n$ as $n \ra \infty$. We
conclude from (3.9) and (3.10):

$$\left|P_0\{\tau_n^\lambda \ is \ even \}- \sum\li_{k=1}^{r(n)}\exp(-\frac{z_k^2}{2\widetilde{m_1}D})
\frac{|I_k|}{\sqrt{2\pi \widetilde{m_1}
D}}\int_{I_k}P_{x+\widetilde{m_1}M} \{\tau_n^\lambda \ is \ even
\}\frac{dx}{|I_k|}\right|\leq \delta+\rho_n^{(4)} \ . \eqno(3.11)$$

Note that, since the random walk $S_n$ is invariant with respect to
shifts, $$P_{x+An}\{\tau_n^\lambda \ is \ even \}=
P_x\{\tau_n^{\lambda+A} \ is \ even\}.$$

Using this fact and taking into account (3.8), we conclude that

$$\int_{I_k}\frac{dx}{|I_k|}P_{x+\widetilde{m_1}M}\{\tau_n^\lambda \ is \ even\}=\int_{I_k}\frac{dx}{|I_k|}
P_x\{\tau_n^{\lambda'} \ is \ even \}+\frac{\rho^{(5)}_n}{\sqrt{n}}
\eqno(3.12)$$ for an appropriate constant $\lambda'$. The integral
in the right hand side of (3.12) is nothing else but the expected
value $\widetilde{\mathbb{E}}$ of the random variable $X$ uniformly
distributed on $I_k$. Let $\chi^{(n)}_x$ be the indicator of the set
$\{\tau_n^{x,\lambda} \ is \ even \}$ in the sample space. Then
$$\int_{I_k}\frac{dx}{|I_k|}P_x\{\tau_n^{\lambda'} \ is \ even \}=\widetilde{\mathbb{E}}\mathbb{E}\chi_X^{(n)}
=\mathbb{E}\widetilde{\mathbb{E}}\chi_X^{(n)}  \ . \eqno(3.13)$$

For given sequences $\{\xi_k\}$, $\{\eta_k\}$,

$$\widetilde{\mathbb{E}}\chi_X^{(n)}=\frac{\eta_1+\eta_2+...+\eta_{\nu_n}+\kappa_n}{\xi_1+\eta_1+\xi_2+\eta_2+...+\xi_{\nu_n}+\eta_{\nu_n}}, \eqno(3.14)$$
where $\nu_n$ and $\kappa_n$ are random variables such that
$\lim\li_{n \ra \infty}\nu_n=\infty$ and $|\kappa_n|<\beta$.

Using the strong law of large numbers, we conclude from (3.14) that
$$\lim\li_{n \ra \infty}\widetilde{\mathbb{E}} \chi_X^n=\frac{\mathbb{E}\eta}{\mathbb{E}\xi+\mathbb{E}\eta} \ , \ a. e. \ .$$

Since $|\widetilde{\mathbb{E}}\chi_X^{(n)}|\leq 1$, the last
equality and (3.13) imply that

$$\lim\li_{n \ra \infty}\int_{I_k}\frac{dx}{|I_k|}P_x\{\tau_n^{\lambda'} \ is \ even \}=\frac{\mathbb{E}\eta}{\mathbb{E}\xi+\mathbb{E}\eta}. \eqno(3.15)$$

From (3.11), (3.12) and (3.15), we derive:

$$\left|P_0\{\tau_n^\lambda \ is \ even \}-\frac{\mathbb{E}\eta}{\mathbb{E}\xi+\mathbb{E}\eta}\right|\leq 2 \delta +\rho_n^{(6)}.$$

This bound implies the first statement of the lemma. The second
statement follows from the first. $\blacksquare$

\

\

Now we have, for our process $\widetilde{X}_t^{\ve,\delta}$,

\

\textbf{Lemma 3.4}

\

\textit{$\lim\limits_{\delta \downarrow 0}\lim\limits_{\varepsilon
\downarrow 0}P\{\widetilde{X}_t^{\varepsilon,\delta}$ finally falls
into the well $\mathcal{E}_1\}=\play{\frac{c_1}{c_1+c_2}}$,}

\textit{$\lim\limits_{\delta \downarrow 0}\lim\limits_{\varepsilon
\downarrow 0}P\{\widetilde{X}_t^{\varepsilon,\delta}$ finally falls
into the well $\mathcal{E}_2\}=\play{\frac{c_2}{c_1+c_2}}$.}

\

(see Fig.7)

\

\textit{Proof:} When the particle collides with the wall $q=-a_1$,
the absolute value of its velocity changes as follows:
$$|p_{new}|=|p_{old}|(1- \ve c_1 - \ve \delta \xi_k).$$

Similarly, when the particle collides with the wall $q=a_2$:
$$|p_{new}|=|p_{old}|(1- \ve c_2 - \ve \delta
\eta_k).$$

Take logarithm of these two equalities:

$$\ln|p_{new}|=\ln|p_{old}|+\ln(1- \ve c_1 - \ve \delta \xi_k),$$

$$\ln|p_{new}|=\ln|p_{old}|+\ln(1- \ve c_2 - \ve \delta \eta_k),$$
respectively. Thus we can view the problem of determining whether
$\widetilde{X}_t^{\ve,\delta}$ enters well $\mathcal{E}_1$ or
$\mathcal{E}_2$ as the following random walk problem. We start from
$\ln|p_0|$, and alternatively jump backward with steplength $-\ln(1-
\ve c_1 - \ve \delta \xi_k)$ (at step $2k-1$) and $-\ln(1- \ve c_2 -
\ve \delta \eta_k)$ (at step $2k$). If the random walk finally jumps
over $\ln|p(O)|$ at an odd step, then $\widetilde{X}_t^{\ve,\delta}$
enters well $\mathcal{E}_1$; otherwise
$\widetilde{X}_t^{\ve,\delta}$ enters well $\mathcal{E}_2$.

Now we further simplify the problem: we start from 0, and
alternatively jump forward with steplength
$U_k^{\ve,\delta}=-\play{\frac{1}{\ve}}\ln(1- \ve c_1 - \ve \delta
\xi_k)$ (at step $2k-1$) and
$V_k^{\ve,\delta}=-\play{\frac{1}{\ve}}\ln(1- \ve c_2 - \ve \delta
\eta_k)$ (at step $2k$). If the random walk finally jumps over
$\play{\frac{1}{\ve}}\ln\play{\left|\frac{p_0}{p(O)}\right|}$ at an
odd step, then $\widetilde{X}_t^{\ve,\delta}$ enters well
$\mathcal{E}_1$; otherwise $\widetilde{X}_t^{\ve,\delta}$ enters
well $\mathcal{E}_2$.

Put $n=\left[\play{\frac{1}{\ve}}\right]$,
$\play{\lambda=\ln\left|\frac{p_0}{p(O)}\right|}$. Taking into
account that $U_k^{\ve,\delta}=(c_1+\delta \xi_k)+O(\ve)$ and
$V_k^{\ve,\delta}=(c_2+\delta \eta_k)+O(\ve)$ as $\ve \da 0$ and
applying Lemma 3.3 we get:

\

$\lim\limits_{\delta \downarrow 0}\lim\limits_{\varepsilon
\downarrow 0}P\{\widetilde{X}_t^{\varepsilon,\delta}$ finally falls
into the well $\mathcal{E}_1\}=\lim\li_{\delta \da 0}\lim\li_{\ve
\da 0}\play{\frac{\mathbb{E}U_k^{\ve, \delta}}{\mathbb{E}U_k^{\ve,
\delta}+\mathbb{E}V_k^{\ve, \delta}}}=\play{\frac{c_1}{c_1+c_2}}.$

\

The second statement of the lemma follows from the
first.$\blacksquare$

\

Let $Y_t=(H_t,K_t)$, $H_0>H(O)$, $K_0=3$, be the process on $\Gamma$
defined inside the edges by formulas given in Corollary (3.1); when
$Y_t$ reaches at time
$t_0=\play{\frac{\sqrt{2}(a_1+a_2)}{c_1+c_2}(H(O)^{-1/2}-H_0^{-1/2})}$
the interior vertex $O$, it goes immediately to $I_1$ or $I_2$ with
probabilities $\play{p_1=\frac{c_1}{c_1+c_2}}$ and
$\play{p_2=\frac{c_2}{c_1+c_2}}$ respectively.

\

Combining Lemma 3.1 and Corollary 3.1 and Lemma 3.4 we have the
following

\

\textbf{Theorem 3.1} \textit{The process
$\widehat{Y}_t^{\varepsilon,\delta}=(\widehat{H}(\widetilde{X}_t^{\varepsilon,\delta}),K(\widetilde{X}_t^{\varepsilon,\delta}))$
on $\Gamma$ converges weakly in uniform topology in the space of
continuous functions $[0,T]\ra \Gamma$ as first $\varepsilon
\downarrow 0$ and then $\delta \downarrow 0$ to the stochastic
process $Y_t=(H(t),K(t))$ on $\Gamma$.}

\

\textbf{Remark:} As we have seen in Section 2, in the case of two
wells, the problem can be regularized by stochastic perturbations of
the initial conditions. Theorem 3.1 shows that the regularization by
perturbations of the dynamics lead to the same limiting slow motion
$Y_t$ on $\Gamma$.

\

Now we consider the case when the model system has more than two
trapping wells (see Fig.8). Then the regularization by perturbation
of the initial conditions, in general, does not work. Let the
coefficients of restitution of the left and right boundaries of the
$i$-th well be $1-\ve c^{(i)}_1 -\ve\delta \xi^{(i)}_k$ for the left
boundary and $1-\ve c^{(i)}_2 -\ve\delta \eta^{(i)}_k$ for the right
boundary. In Fig.8 $i \in \{1,2,3,4,5\}$. A random variable with the
subscript $k$ is used when the particle hits a wall the $k$-th time;
$c^{(i)}_1$, $c^{(i)}_2$ are positive constants;
$\{\xi^{(i)}_k\}_{k=1}^\infty$, $\{\eta^{(i)}_k\}_{k=1}^\infty$ are
sequences of i.i.d random variables. We write $\xi^{(i)}$,
$\eta^{(i)}$ to denote independent random variables that have the
same distributions as any of the random variables from the sequences
$\{\xi^{(i)}_k\}_{k=1}^{\infty}$, $\{\eta^{(i)}_k\}_{k=1}^{\infty}$.
We assume that
$P\{\alpha<\xi^{(i)}<\beta\}=P\{\alpha<\eta^{(i)}<\beta\}=1$ for
$0<\alpha<\beta<\infty$.

The resulting process is denoted by
$X^{\varepsilon,\delta}_t=(q^{\varepsilon,\delta}(t),p^{\varepsilon,\delta}(t))$.
We rescale time and put
$\widetilde{X}^{\varepsilon,\delta}_t=X^{\varepsilon,\delta}_{t/\varepsilon}$.
Let
$Y_t^{\varepsilon,\delta}=(H(\widetilde{X}_t^{\varepsilon,\delta}),K(\widetilde{X}_t^{\varepsilon,\delta}))$
be the projection of $\widetilde{X}_t^{\ve, \delta}$ onto the graph
$\Gamma$ (see Fig.8). Put
$H^{\varepsilon,\delta}(t)=H(\widetilde{X}_t^{\varepsilon,\delta})$.
For fixed $\varepsilon >0$, this is a stochastic process with jumps
and each sample path is a step function which is right continuous
and has limit from the left. Using the same construction as in the
end of Section 2, we obtain from $H^{\ve, \delta}(t)$ a continuous
piecewise linear process $\widehat{H}^{\varepsilon, \delta}(t)$. We
also put
$\widehat{Y}_t^{\varepsilon,\delta}=(\widehat{H}(\widetilde{X}_t^{\varepsilon,\delta}),K(\widetilde{X}_t^{\varepsilon,\delta}))$.

\begin{figure}
\centering
\includegraphics[height=7cm, width=15cm , bb=44 58 513 235]{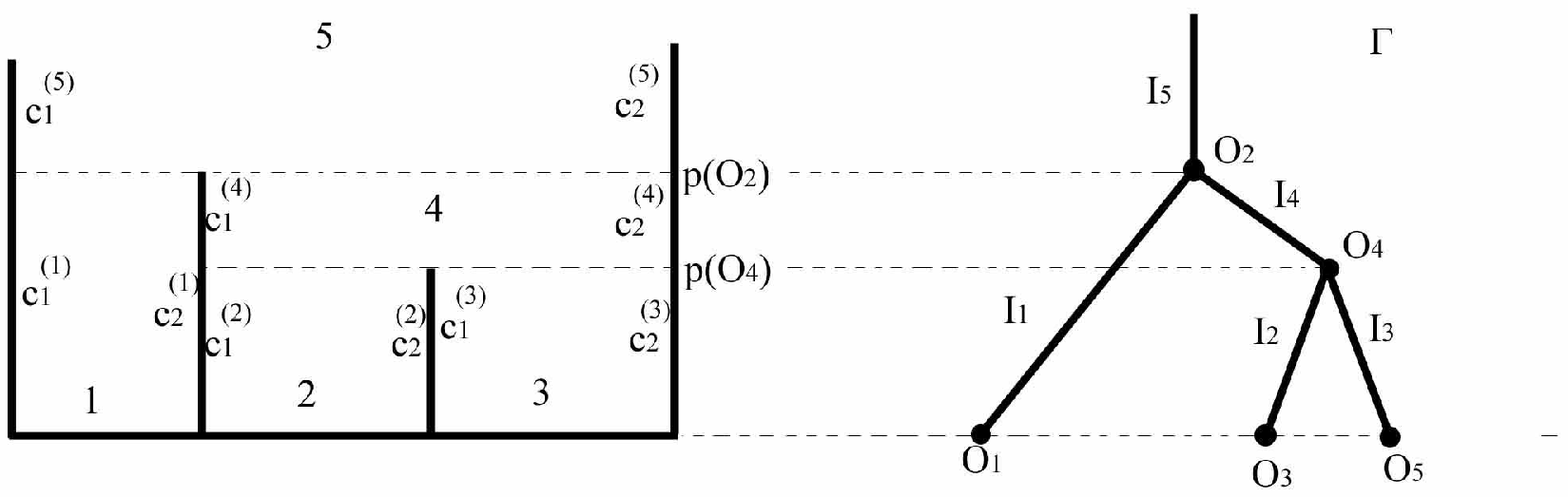}
\caption{The case when the system has more than two trapping wells}
\end{figure}

The averaging procedure is the same as before. One can see that
within each edge $I_i$, as first $\ve \da 0$ and then $\delta \da
0$, the process $H^{\ve,\delta}(t)$ converges to a deterministic
motion $H(t)$ which satisfies the differential equation
$$\play{\frac{dH}{dt}=-2\frac{c_1^{(i)}+c_2^{(i)}}{T_i(H)}H}, \eqno(3.16)$$
where $i \in \{1,2,3,4,5\}$. $T_i(H)$ is the period of motion in the
phase space for the particle in the $i$-th well with energy $H$:
$T_i(H)=\play{\frac{\sqrt{2}m_i}{\sqrt{H}}}$ where $m_i$ is the
width of the $i$-th well.

The branching probabilities for the limiting motion $Y(t)=(H(t),
K(t))$ at each interior vertex $O_l$ are given by
$\play{p_1^{(l)}=\frac{c^{(i)}_1}{c^{(i)}_1+c^{(i)}_2}}$ and
$\play{p_2^{(l)}=\frac{c^{(i)}_2}{c^{(i)}_1+c^{(i)}_2}}$. Here
$i=i(l)$ is the number of the edge that is above $O_l$; $p_1^{(l)}$,
$p_2^{(l)}$ are probabilities that the particle enters the left edge
or right edge that is below $O_l$. The motion inside the edges is
described as in Corollary (3.1). Since each time when the particle
hits the boundary, the perturbation of the dynamics is given by
independent random variables, the branching at each interior vertex
is also independent of each other. Therefore we finally have the
following theorem:

\

\

\textbf{Theorem 3.2} \textit{The process
$\widehat{Y}_t^{\varepsilon,\delta}=(\widehat{H}(\widetilde{X}_t^{\varepsilon,\delta}),K(\widetilde{X}_t^{\varepsilon,\delta}))$
on $\Gamma$ converges weakly in uniform topology in the space of
continuous functions $[0,T]\ra \Gamma$ as first $\varepsilon
\downarrow 0$ and then $\delta \downarrow 0$ to the continuous
stochastic process $Y_t=(H(t),K(t))$ on $\Gamma$. The process $Y_t$,
$Y_0=Y(q_0,p_0)$, inside the edges is described by (3.16). When the
process $Y_t$ along an edge $I_i$ attached to an interior vertex
$O_l$ and situated above $O_l$ reaches $O_l$, it instantaneously
leaves $O_l$ and enters edges of $\Gamma$ attached to $O_l$ and
situated below $O_l$ on the left and on the right with probabilities
$\play{p_1^{(l)}=\frac{c^{(i)}_1}{c^{(i)}_1+c^{(i)}_2}}$,
$\play{p_2^{(l)}=\frac{c^{(i)}_2}{c^{(i)}_1+c^{(i)}_2}}$,
respectively. The branching at each interior vertex $O_l$ is
independent of each other. }

\

\textbf{Remark:} Let us suppose that $\delta=1$: the coefficients of
restitution become $1-\ve c^{(i)}_1-\ve \xi^{(i)}_k$ and $1-\ve
c^{(i)}_2-\ve \eta^{(i)}_k$, respectively. Then the process
$\widetilde{Y}_t^{\ve,1}$ converges weakly as $\ve \da 0$ to the
process $\widetilde{Y}_t$ on $\Gamma$ described above with the
replacement of $c_1^{(i)}$ and $c_2^{(i)}$ by
$\widetilde{c}_1^{(i)}=c_1^{(i)}+\mathbb{E}\xi^{(i)}$ and by
$\widetilde{c}_2^{(i)}=c_2^{(i)}+\mathbb{E}\eta^{(i)}$,
respectively.

\

\

\section{The case of general potential}

The problem mentioned in the beginning of the introduction is
considered in this section. Let the motion of a particle inside
$[a_1,a_2]$ be governed by the equation
$\ddot{q}_t^\ve=-F'(q_t^\ve)$, $q_0^\ve=q_0$, $\dot{q}_0^\ve=p_0$.
When the particle hits the walls, it is reflected with a loss of
energy. Let $H(p,q)=\play{\frac{p^2}{2}+F(q)}$ be the energy of the
particle, $p=\dot{q}$. The coefficients of restitution at $a_i$ are
equal to $1-\ve c_i(H)$, $0<\ve<<1$, $i\in \{1,2\}$, $c_i(H)$ are
positive and smooth. As was explained in the introduction, we can
restrict ourselves to the case when $F(q)$ has just one maximum
inside $[a_1,a_2]$, say, at $a_0\in (a_1,a_2)$ (see Fig.9). Let
$\Gamma$ be the graph corresponding to the Hamiltonian $H(p,q)$,
$\sqcap=\{(p,q)\in \mathbb{R}^2: a_1\leq q \leq a_2\}$, and $Y:
\sqcap \ra \Gamma$ be the projection of $\sqcap$ on $\Gamma$:
$Y(p,q)=(H(p,q),K(p,q))$.

\begin{figure}
\centering
\includegraphics[height=6cm, width=10cm , bb=55 45 512 283]{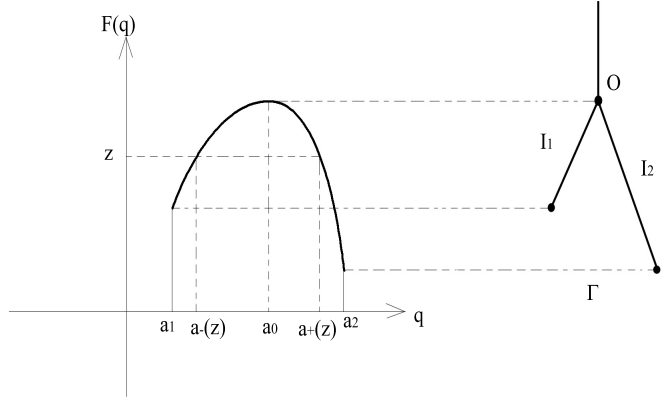}
\caption{The case of general potential}
\end{figure}

Denote by $X_t^\ve=(q_t^\ve,p_t^\ve)$ nearly-elastic trajectory
starting at $(q_0,p_0)$. This motion has a fast component, which is
actually the motion along the elastic trajectories, and the slow
component $Y_t^\ve=Y(X_t^\ve)$. Let
$\widetilde{Y}_t^\ve=Y^\ve_{t/\ve}$ and let
$\widehat{Y}^\ve_{t/\ve}$ be the continuous piecewise linear
approximation of $Y^{\ve}_{t/\ve}$ as we considered earlier. As in
the model problem, one can see that $\lim\li_{\ve \da
0}\widehat{Y}_t^\ve$ (and $\lim\li_{\ve \da 0}\widetilde{Y}_t^\ve$)
does not exist. Since we have just two wells, one can use
perturbations of the initial conditions to regularize the problem.
Let $\xi_\delta$ be a two-dimensional random variable with a
continuous density which is positive for $|x|<\delta$ and is equal
to zero if $|x|>\delta$. Denote by $\widehat{Y}_t^{\ve,\delta}$ the
continuous modification introduced above when the initial condition
$x_0=(q_0,p_0)$ is replaced by $x_0+\xi_\delta$;
$\widehat{Y}^{\ve,\delta}_t$ is a stochastic process since
$\widehat{Y}_0^{\ve,\delta}$ is random.

Let us introduce now a stochastic process $Y_t=(H_t,K_t)$ on
$\Gamma$. First, let $T_k(H)$, $k \in \{1,2,3\}$, be the period of
oscillations of the elastic motion $q_t^0$ with the initial
conditions in $Y^{-1}(H,K)$. It is easy to check that

$$T_3(H)=2\int_{a_1}^{a_2}\frac{dq}{\sqrt{2(H-F(q))}} \ ,$$ $$T_1(H)=2\int_{a_1}^{a_-(H)}\frac{dq}{\sqrt{2(H-F(q))}} \
,$$ $$T_2(H)=2\int_{a_+(H)}^{a_2}\frac{dq}{\sqrt{2(H-F(q))}} \ .$$

Here $a_{\pm}(z)$ are roots of the equation $F(a_{\pm}(z))=z$,
$a_-(z)<a_+(z)$.

Let $Y_t=(H_t,K_t)$ be the stochastic process on $\Gamma$ such that
$H_0=H(q_0,p_0)>H(O)$, $K_0=3$; $H_t$ is a deterministic motion
inside each edge:
$$\dot{H}_t=-2\play{\frac{c_k^{(1)}(H_t)(H_t-F(a_k^{(1)}))+c_k^{(2)}(H_t)(H_t-F(a_k^{(2)}))}{T_k(H_t)}}$$
inside $I_k$, $k \in \{1,2,3\}$, where $c_k^{(j)}(H)=c_j(H)$,
$a_k^{(j)}=a_j$, for $k=3$ and $j=1,2$; $c_k^{(j)}(H)=c_k(H)$,
$a_k^{(j)}=a_k$ for $k=1,2$ and $j=1,2$. The trajectory $Y_t$ hits
$O$ in a finite time
$t_0=\play{\int_{H(O)}^{H(q_0,p_0)}\frac{T_3(z)dz}{2[c_1(z)(z-F(a_1))+c_2(z)(z-F(a_2))]}}$
since $T_3(H)\sim \ln (H-H(O))$ as $H \ra H(O)$. After hitting
vertex $O$, $Y_t$ leaves $O$ immediately and goes to $I_1$ or $I_2$
with probabilities
$p_1=\play{\frac{c_1(H(O))}{c_1(H(O))+c_2(H(O))}}$ and
$p_2=\play{\frac{c_2(H(O))}{c_1(H(O))+c_2(H(O))}}$. These conditions
define $Y_t$ in a unique way.

\

\textbf{Theorem 4.1} \textit{Under mentioned above conditions,
process $\widehat{Y}_t^{\ve,\delta}$ converge weakly in the space
$C_{0T}$ of continuous functions $[0,T]\ra \Gamma$, $T<\infty$, as
first $\ve \da 0$ and then $\delta \da 0$ to the process $Y_t$.}

\

The proof of this Theorem is similar to the proof of Theorem 2.1 and
we omit it.

\

One can regularize the problem by introducing random perturbations
not of the initial conditions, but of the dynamics, as we did in
Section 3. Similar to Theorem 3.1, one can prove that the process
$Y_t$ introduced above again serves as the limiting slow motion.
These results allow to say that stochasticity is an intrinsic
property of the long-time behavior of nearly-elastic systems: the
limiting slow motion is the same stochastic process for different
regularizations.

\

\section{A two dimensional model problem}

Consider a standard billiard model in which a point mass moves
freely inside a convex, bounded, and simply connected region $G$ in
the plane with smooth boundary $\partial G$ (Fig.10). The orbits of
such motion consist of straight line segments inside $G$ joined at
boundary points according to the rule that the angle of incidence
equals the angle of reflection. Speed is a constant of motion. The
phase space $\bigwedge$ of this system is conveniently described as
the set of all tangent vectors of fixed length (say, unit length)
supported at points of the interior of $G$ together with vectors at
boundary points pointing inward. We parameterize the phase space
$\bigwedge$ by the coordinate $(g,\varphi)$, where $g\in G\cup
\partial G$ is the position of the particle within $G \cup \partial G$ and
$\varphi\in [0,2\pi)$ is the angle between the positive $Y$-axis and
the velocity vector of the particle. The phase space of the motion
has a global section $\Pi$ which consists of all inward-pointing
unit vectors on $\partial G$. Topologically $\Pi$ is a cylinder
parameterized by the cyclic length parameter $s\in (0,L]$ along
$\partial G$ (here $L$ is the length of $\partial G$) and the angle
$\theta \in [0,\pi]$ between the velocity and the positive
tangential direction. The Poincar\'{e} map $f: \Pi \rightarrow \Pi$
is usually called \textit{the billiard map} and can be described by
$f(s,\theta)=(S(s,\theta), \Theta(s,\theta))$. Notice that when
$\partial G$ is smooth $f$ is also smooth. For a further reference
about the billiard model we refer to [9], page 339.

A classical result in such a billiard model is that $f$ preserves
the volume element $\sin \theta dsd\theta$ on $\Pi$. This volume
element $\sin \theta dsd\theta$ is usually called \textit{the
Liouville measure}; $m(s,\theta)=\sin\theta$ is the density of the
Liouville measure.

\begin{figure}
\centering
\includegraphics[height=9cm, width=15cm , bb=68 21 542 281]{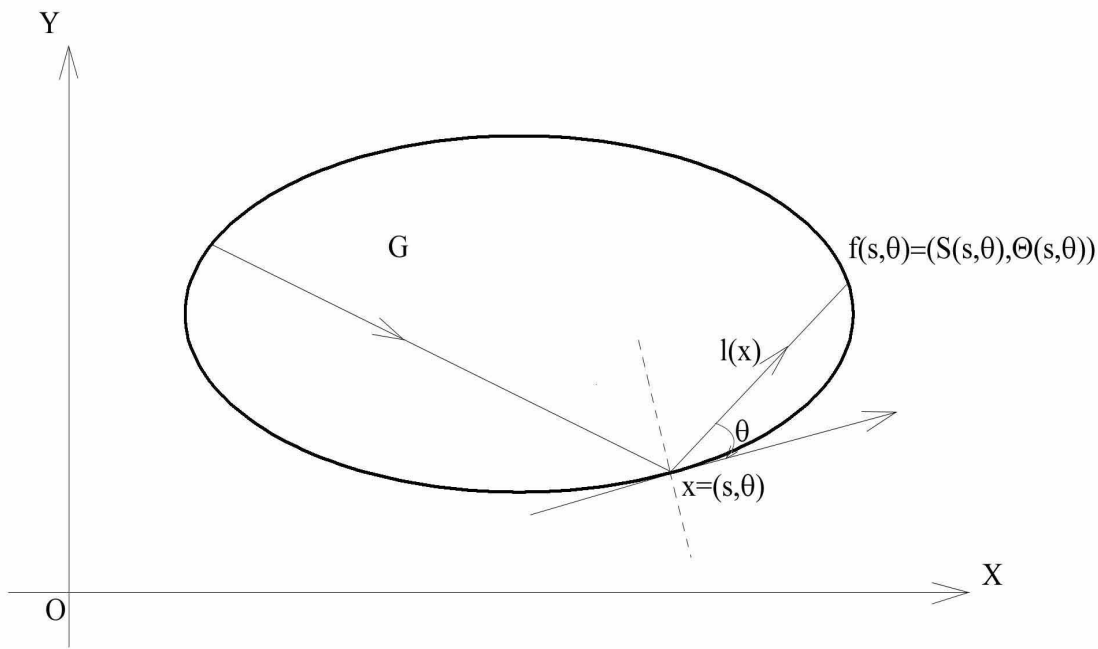}
\caption{The billiard model}
\end{figure}

We now assume that the particle loses a small amount of energy each
time when it has a collision with $\partial G$. Let $c(s,\theta)$ be
a smooth function on $\Pi$, bounded from above and below by some
fixed positive constants. We assume for brevity that the energy $H$
of the particle after the collision
 with $\partial G$ at a point $x\in \Pi$ diminishes to $\max\{H-\varepsilon c(x),0\}$
 (rather than to $H(1-\ve c(x))$ as before); here $\ve$ is a small positive parameter. The
direction of the velocity vector remains the same as for the elastic
system (see Fig.10).

We assume that there is another wall of height $H(O)$ separating the
region $G$ (see Fig.11). In a finite time of order
$\play{\frac{const}{\ve}}$, $\ve \da 0$, each trajectory of the
nearly-elastic motion, which had at time zero energy $H_0>H(O)$,
enters one of the wells 1 or 2 shown in Fig.11 and continues the
motion there (we assume that the loss of energy on the wall
separating the region $G$ is defined in a similar way). Which of the
wells is entered, in general, depends on $\ve$ in a very sensitive
way: as $\ve \da 0$ the trajectory alternatively enters well 1 or 2.
As in the previous sections, the nearly-elastic motion has a fast
and a slow components. In the case of one degree of freedom, the
fast motion has just one normalized invariant measure, and a unique,
independent of the initial conditions, limiting slow motion inside
the edges exists without any regularization. This is not the case
now. Besides the Liouville measure, other invariant measures on the
Poincar\'{e} section $\Pi$ exist; for instance, there are periodic
points on $\partial G$. On the other hand, it is natural to assume
that the reflections on $\partial G$ are undergoing small random
perturbations of intensity $\delta << 1$ so that the resulting
motion is a stochastic process depending on two small parameters
$\ve$ and $\delta$. We choose a class of random perturbations in
such a way that the perturbed system induces on $\Pi$ a unique
invariant measure, namely, the Liouville measure (compare with [5]).
It turns out that the perturbations in the fast motion which we
introduce provide a regularization of the slow component not just
inside of the edges of the corresponding graph but also near the
vertices, so that the weak limit of the slow component of the
regularized motion (after an appropriate time change) exists on the
whole graph as first $\ve$ and then $\delta$ tends to zero.
Probably, the class of permissible random perturbations leading to
the same limiting slow motion can be essentially extended. But one
should keep in mind that the existence of a rich set of invariant
measures for the billiard map shows that the limiting slow motion
calculated in this section has, in a sense, a restricted
universality: other regularizations of the fast motion may lead to
different slow motions.

We regularize this problem by considering the following scheme of
small stochastic perturbation of the system. Consider an
one-dimensional diffusion process $I_{t}^{\theta}$ on $[0,\pi]$
starting from point $\theta$. The diffusion process is governed by a
second order elliptic operator $\mathcal{L}$ which is self-adjoint
with respect to the Liouville measure. Such an operator can be
written as
$\mathcal{L}u=\displaystyle{\frac{1}{2\sin{\theta}}\frac{d}{d\theta}(a(\theta)\frac{du}{d\theta})}$,
where $a$ is a smooth positive function of $\theta$. We assume that
the diffusion process $I_t^\theta$ has instantaneous reflection at
the boundary points $\theta=0, \pi$. Notice that there is a
singularity of our operator $\mathcal{L}$ at the boundary points
$\theta=0, \pi$. Inspite of this singularity the boundary points
$\theta= 0, \pi$ are accessible and therefore we need boundary
conditions: we choose instantaneous reflection at these points.

Choose a small constant $\delta
>0$. Let the position of the particle after one collision be $x=(s,\theta)\in \Pi$ with an energy $H>0$. Then the
next position on $\Pi$ will be not exactly $x$ but
$Z_\delta^x=(s,I_\delta^\theta)$, and the energy will diminish to
$H-\ve c(x)$ (see Fig.11). And the subsequent motion starts from
this new point $Z_\delta^x$. The same kind of perturbation is
repeated again and again for each collision, independent of other
collisions. This scheme of perturbation models the fact that each
time when a collision happens, the particle changes its direction of
motion a little bit, due to the random perturbations.

After adding the small stochastic perturbation, the evolution of
energy and the motion of the particle is described by a continuous
time Markov process $N_t^{\varepsilon, \delta}=(H_t^{\varepsilon,
\delta}, \Lambda_t^{\ve,\delta})$ on $[0,+\infty)\times \bigwedge$.
The energy $H_t^{\ve,\delta}$ changes just on the boundary: if
$\Lambda_t^{\ve,\delta}=x=(s,\theta)\in \Pi$, the energy
instantaneously decreases on $\ve c(x)$; $H_t^{\ve,\delta}$ is right
continuous. The process $\Lambda_t^{\ve,\delta}=(g_t^{\ve,\delta},
\varphi_t^{\ve,\delta})$ has a piecewise linear first component: it
moves uniformly along the chords connecting successive points where
$g_t^{\ve,\delta}$ hits $\partial G$ with the constant speed
$\sqrt{2H_t^{\ve,\delta}}$. Let $M_t^{\ve,\delta}=(H_t^{\ve,\delta},
g_t^{\ve,\delta})$  (a sample trajectory of $M_t^{\ve,\delta}$ is
the thin line in Fig.11). The trajectory of $M_t^{\ve,\delta}$ is
right continuous and has limit from the left. We define a continuous
modification
$\widehat{M}_t^{\ve,\delta}=(\widehat{H}_t^{\ve,\delta},
g_t^{\ve,\delta})$ of $M_t^{\ve,\delta}$ by taking
$\widehat{H}_t^{\ve,\delta}$ as a piecewise linear modification of
$H_t^{\ve,\delta}$ (a sample trajectory of
$\widehat{M}_t^{\ve,\delta}$ is the bold line in Fig.11). Let
$\widehat{N}_t^{\ve,\delta}=(\widehat{H}_t^{\ve,\delta},
\Lambda_t^{\ve,\delta})$.

\begin{figure}
\centering
\includegraphics[height=9cm, width=15cm , bb=76 2 542 260]{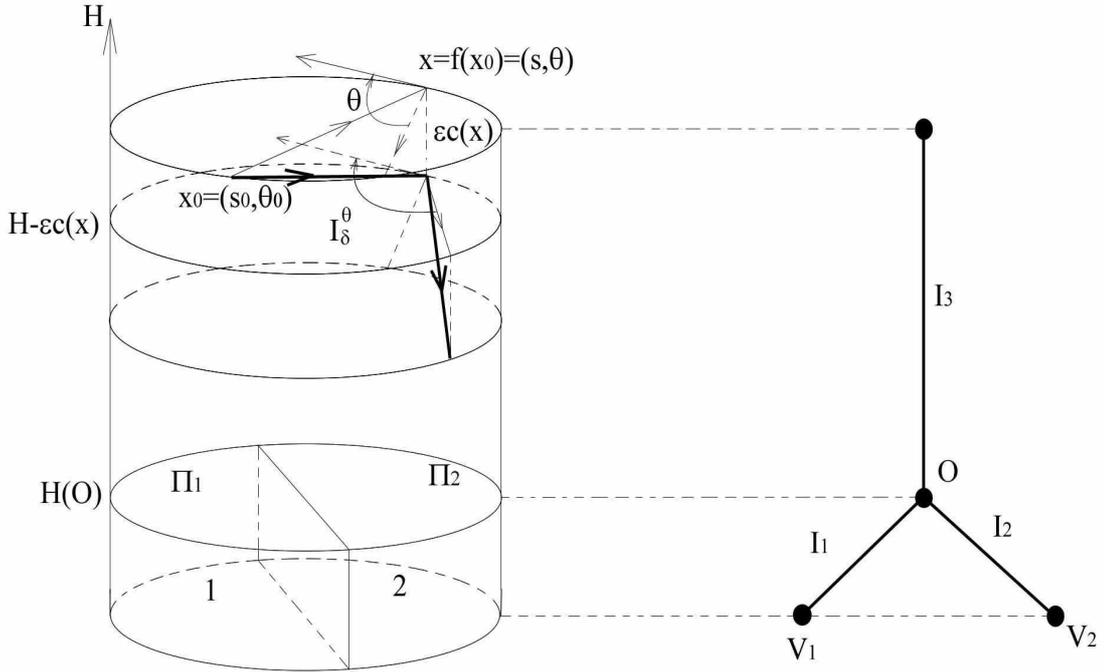}
\caption{The billiard model with energy loss and perturbation}
\end{figure}

Consider also a jumping Markov chain $\{X_n^\delta\}_{n\geq 1}$ on
$\Pi$, defined as follows: from a point
$X_n^\delta=(s_n^\delta,\theta_n^\delta)\in \Pi$ the chain goes in
one time unit to
$X_{n+1}^\delta=(s_{n+1}^\delta,\theta_{n+1}^\delta)=(S(s_n^\delta,\theta_n^\delta),I_\delta^{\Theta(s_n^\delta,\theta_n^\delta)})\in
\Pi$.

Suppose the particle of unit mass starts its motion from a point
$x=(s,\theta)\in \Pi$ and the energy of the particle is $H$. Let
$T(x,H)$ be the time that the particle starting from $x$ first
reaches $\partial G$. Let $l(x)$ be the chord starting from $x\in
\Pi$ (see Fig.10). Let $L(x)$ be the length of $l(x)$. We have
$T(x,H)=\play{\frac{L(x)}{\sqrt{2H}}}$.

We identify points of each well having the same energy. The set
obtained after such an identification is a graph $\Gamma$; the
metric on $\Gamma$ is defined by the distance along the edges of
$\Gamma$. This graph has one interior vertex $O$ corresponding to
the energy level $H(O)$ connected with three edges: $I_3$
corresponding to trajectories with energy greater than $H(O)$, $I_1$
corresponding to trajectories with energy less than $H(O)$ situated
within well 1 (to the left of $O$), and $I_2$ corresponding to
trajectories with energy less than $H(O)$ situated within well 2 (to
the right of $O$).

We consider the process $Y_t^{\varepsilon,
\delta}=(H_t^{\varepsilon, \delta},K_t^{\varepsilon, \delta})$ on
the tree $\Gamma$, where $K_t^{\varepsilon, \delta}$ is defined as
follows: $K_t^{\varepsilon, \delta}=3$ if the energy of the particle
is greater than $H(O)$, $K_t^{\varepsilon, \delta}=1$ if the
particle has energy less than $H(O)$, while
$N_t^{\varepsilon,\delta}$ is in the potential well 1, and
$K_t^{\varepsilon, \delta}=2$ if the particle has energy less than
$H(O)$, while $N_t^{\varepsilon,\delta}$ is in the potential well 2.
The process $Y_t^{\ve,\delta}$ is the slow component of
$N_t^{\ve,\delta}$ as $\ve \da 0$. We show that the rescaled process
$Y_{t/\varepsilon}^{\varepsilon, \delta}$ converges, as first
$\varepsilon \downarrow 0$ then $\delta \downarrow 0$, to a
stochastic process $Y_t$ on the tree $\Gamma$. Such a stochastic
process has deterministic motion within each edge of the tree
$\Gamma$ and has stochasticity only at the interior vertex $O$.

\

Let us first state some auxiliary results:

\

\textbf{Lemma 5.1} \textit{The diffusion process $I_t^{\theta}$ on
$[0,\pi]$ governed by the operator
$\mathcal{L}u=\displaystyle{\frac{1}{2\sin{\theta}}\frac{d}{d\theta}(a(\theta)\frac{du}{d\theta})}$
with instantaneous reflection at boundary points $\theta=0, \pi$,
has $\play{\frac{1}{2}\sin\theta}$ as its unique invariant density.}

\

Suppose small $r>0$, let $\Pi(r)=\{x=(s,\theta)\in \Pi;
r\leq\theta\leq \pi-r \}$.

\

\textbf{Lemma 5.2} \textit{For any $\delta > 0$, for any small
enough $r>0$, the
 Markov chain $\{X_n^\delta\}_{n \geq 1}$ on
$\Pi(r)$ satisfies the Doeblin condition, and has only one ergodic
component on $\Pi(r)$. The invariant measure for the process
$\{X_n^{\delta}\}_{n\geq 1}$ on $\Pi$ has a density
$\play{d(x)=\frac{1}{2L}m(x)=\frac{1}{2L}\sin\theta}$.}

\

\textbf{Lemma 5.3} \textit{Within edge $I_3$ of the graph $\Gamma$,
$\lim\limits_{\delta \downarrow 0}\lim\limits_{\varepsilon
\downarrow 0}H^{\varepsilon, \delta}_{t/\varepsilon}=H_t$ in
probability, where $H_t$ satisfies the following differential
equation:}

$$\displaystyle{\frac{dH_t}{dt}=-\frac{\displaystyle{\iint_{\Pi}
c(x)m(x)dx}}{\displaystyle{\iint_{\Pi}T(x,H_t)m(x)dx}}} \ , \
H_0=H_0^{\ve,\delta} \eqno(5.1)$$

\

In exactly the same way as this lemma one can have similar results
for limiting deterministic motion within edge $I_1$ and $I_2$. We
omit details here.

\

\

\textbf{Lemma 5.4} \textit{Denote the area of the region bounded by
the convex curve in Fig.10 by A. Then $$\iint_\Pi L(x)m(x)dx=2\pi A
\,$$ where $m(x)$ is the density of the Liouville measure on $\Pi$
and $L(x)$ is the length of the chord $l(x)$.}

\

Using Lemma 5.4,  equation (5.1) can be simplified as
$$\frac{dH_t}{dt}=-\sqrt{2H_t}\frac{1}{2 \pi A}\iint_\Pi c(x)m(x)dx. \eqno(5.2)$$

And (5.2) has an explicit solution

$$H_t=\left(\sqrt{H_0}-\frac{t}{2\sqrt{2}\pi A}\iint_\Pi c(x)m(x)dx\right)^2 \ ,
\ 0 \leq t \leq \frac{2\pi A
(\sqrt{2H_0}-\sqrt{2H(O)})}{\play{\iint_\Pi c(x)m(x)dx}} \ .
\eqno(5.3)$$

\

\textbf{Lemma 5.5} \textit{For any $\ve>0$, any $\delta>0$ the
invariant measure for the process $\Lambda_t^{\ve,\delta}$ is
proportional to the Lebesgue measure on $\bigwedge$.}

\

\textbf{Lemma 5.6} \textit{For any $\delta > 0$,}

\

\textit{$\lim\li_{\ve \da 0}P\{N_t^{\ve,\delta}$ finally falls into
well 1$\}= \lim\li_{\ve \da 0}P\{\widehat{N}_t^{\ve,\delta} $finally
falls into well 1$\}$,}

\textit{$\lim\li_{\ve \da 0}P\{N_t^{\ve,\delta}$ finally falls into
well 2$\}= \lim\li_{\ve \da 0}P\{\widehat{N}_t^{\ve,\delta} $finally
falls into well 2$\}$. }

\

At the beginning of this section we have already defined the process
$Y_t^{\varepsilon, \delta}=(H_t^{\varepsilon,
\delta},K_t^{\varepsilon, \delta})$ on the tree $\Gamma$. Let
$\Pi_i, \ i \in \{1,2\}$ be the set of those points on $\Pi$ which
correspond to well $i \ , \ i \in \{1,2\}$, respectively (see
Fig.11, but notice that it only shows the $s$ coordinate in the
space $\Pi$).

\

\textbf{Lemma 5.7}

$\lim\li_{\delta \da 0}\lim\li_{\ve \da 0} P\{N_t^{\ve,\delta}$
\textit{finally falls into well 1}
$\}=\play{\frac{\play{\iint_{\Pi_1}c(s,\theta)m(s,\theta)dsd\theta}}
{\play{\iint_{\Pi}c(s,\theta)m(s,\theta)dsd\theta}}} \ ,$

$\lim\li_{\delta \da 0}\lim\li_{\ve \da 0} P\{N_t^{\ve,\delta}$
\textit{finally falls into well 2}
$\}=\play{\frac{\play{\iint_{\Pi_2}c(s,\theta)m(s,\theta)dsd\theta}}
{\play{\iint_{\Pi}c(s,\theta)m(s,\theta)dsd\theta}}} \ .$

\

\

Now similarly as in Section 2 and Section 3,  we consider a
piecewise linear modification of the process
$Y_t^{\varepsilon,\delta}$, and denote such a process by
$\widehat{Y}_t^{\varepsilon,\delta}$. We have
$\widehat{Y}_t^{\ve,\delta}=(\widehat{H}_t^{\ve,\delta},K_t^{\ve,\delta})$.
For all $ \varepsilon
>0$, for all $0<t<T$, for some uniform constant $C>0$,
$P(|\widehat{Y}_t^{\varepsilon,\delta}-Y_t^{\varepsilon,\delta}|<C\varepsilon)=1$.
We define a stochastic process $Y_t$ on $\Gamma$ as follows: it has
deterministic motion within each edge of the tree $\Gamma$ and only
has stochasticity at the interior vertex $O$. Let $Y_t=(H(t),K(t))$
where $H(t)$ is the energy of the particle and $K(t)$ is the
 number of the edge of the graph $\Gamma$. The
deterministic motion within each edge of the tree $\Gamma$ is
defined as follows: on edge 3 the first component $H(t)$ of $Y_t$
satisfies
$$H_t=\left(\sqrt{H_0}-\frac{t}{2\sqrt{2}\pi A}\iint_\Pi c(x)m(x)dx\right)^2 \ , \ 0 \leq t \leq t_0 \ , \eqno(5.4)$$
and on edge 1 and 2 the differential equation becomes

$$H_t=\left(\sqrt{H(O)}-\frac{t-t_0}{2\sqrt{2}\pi A_1}\iint_{\Pi_1} c(x)m(x)dx\right)^2 \ , \ t_0 \leq t \leq t_1 \ , \eqno(5.5)$$

$$H_t=\left(\sqrt{H(O)}-\frac{t-t_0}{2\sqrt{2}\pi A_2}\iint_{\Pi_2} c(x)m(x)dx\right)^2 \ , \ t_0 \leq t \leq t_2 \ , \eqno(5.6)$$
where $\Pi_1$ and $\Pi_2$ are defined as before, and $A_1$, $A_2$
are the areas of the domains in $G$ corresponding to well 1 and 2.
The process $Y_t$, starting from a point $Y_0=(H_0, 3)$, will reach
the interior vertex $O$ in time $t_0=\play{\frac{2\pi
A(\sqrt{2H_0}-\sqrt{2H(O)})}{\play{\iint_\Pi c(x)m(x)dx}}}<\infty$
and will instantaneously leave $O$ and enter one of the edges 1 or 2
with probabilities
$p_1=\displaystyle{\frac{\displaystyle{\iint_{\Pi_1}c(x)m(x)dx}}{\displaystyle{\iint_{\Pi}c(x)m(x)dx}}}$
and
$p_2=\displaystyle{\frac{\displaystyle{\iint_{\Pi_2}c(x)m(x)dx}}{\displaystyle{\iint_{\Pi}c(x)m(x)dx}}}$
respectively. After $Y_t$ enters edge 1 or 2 it will move
deterministically according to the given differential equation (5.5)
or (5.6) till time $\play{t_i=t_0+\frac{2\pi
A_i\sqrt{2H(O)}}{\play{\iint_{\Pi_i}c(x)m(x)dx}}}$, $i \in \{1,2\}$,
respectively, when it hits energy level zero.

\

These lemmas imply the following:

\

\textbf{Theorem 5.1} \textit{The process
$\widehat{Y}_{t/\varepsilon}^{\varepsilon,\delta}$ converges weakly,
as first $\varepsilon \downarrow 0$ then $\delta \downarrow 0$, to
$Y_t$. }

\

Let us prove now these lemmas.

\

\textit{Proof of Lemma 5.1:} Let $\tau$ be the first time for
process $I_t^\theta$ to exit from the interval $[0,\pi]$. By
considering $u(\theta)=\mathbb{E}_\theta\tau$ as the solution of the
corresponding Dirichlet problem, one can check that $u(\theta)$ is
finite so that the boundary points $\theta=0, \pi$ are accessible.

 Now we consider another diffusion process
$\widetilde{I}_t^{\theta}$ which starts from point $\theta$ and is
governed by the operator
$\displaystyle{\widetilde{\mathcal{L}}f=\frac{1}{2}\frac{d}{d\theta}(a(\theta)\frac{df}{d\theta})}$
with instantaneous reflection at the boundary points $\theta=0,
\pi$. The invariant denisty for the process
$\widetilde{I}_t^{\theta}$ is the uniform distribution with the
density $\play{\frac{1}{\pi}}$.

Our process $I_t^{\theta}$ can be obtained from
$\widetilde{I}_t^{\theta}$ by taking a time change $dt=\sin\theta
d\widetilde{t}$. Therefore invariant measure for $I_t^{\theta}$ is
$\play{\frac{1}{2}}\sin \theta d\theta$ (see, [8], Chapter 5 for a
reference concerning the random time change). $\blacksquare$

\

\textit{Proof of Lemma 5.2:} For any small $r>0$, any $\delta>0$,
any $\theta_1, \theta_2\in [r,\pi-r]$, we have $p_I(\delta;
\theta_1, \theta_2)\geq a(r)>0$ for a positive constant $a(r)>0$
depending on $r$, where $p_I(\cdot; \cdot, \cdot)$ is the transition
density for the process $I_t^\theta$.

Let $r>0$. Let $\Pi(r)=\{x=(s,\theta)\in \Pi; r\leq \theta \leq
\pi-r\}$. Since the function $\Theta(x)\equiv\Theta(s,\theta)$ is
continuous for $x\in \Pi(r)$ and the set $\Pi(r)$ is compact,
$\Theta(x)$ must attain its maximum and minimum on $\Pi(r)$. It is
easy to see that there exist $\lambda(r)>0$ such that $\min(r,
\pi-\max\li_{x\in \Pi(r)}\Theta(x), \min\li_{x\in
\Pi(r)}\Theta(x))=\lambda(r)>0$ (otherwise the particle must start
from the tangential direction and keep on moving along the boundary,
which contradicts the fact that $\theta\in [r,\pi-r]$). This implies
that $f(\Pi(r))\subseteq \Pi(\lambda(r))$ with $0<\lambda(r)\leq r$.

We check the following form of Doeblin condition: For any $x\in
\Pi(r)$, there exist an integer $n_0\geq 1$ and a positive $\ve>0$
such that for any Borel measurable function $A$ in
$\mathcal{B}(\Pi(r))$, $\nu(A)\leq \ve$, we have $p^{(n_0)}(x,A)\leq
1-\ve$. Here $\nu(\cdot)$ is the standard Lebesgue measure on $\Pi$.

Let $x=(s_1,\theta_1)\in \Pi(r)$. Let
$(s_2,\theta_2)=(S(s_1,\theta_1), \Theta(s_1,\theta_1))\in
\Pi(\lambda(r))$. Let a point $(s_3,\theta_3)\in \Pi(r)$, $s_3\neq
s_2$. Let $\zeta(s_2, s_3)\in (0, \pi)$ be the unique angle such
that $S(s_2,\zeta(s_2, s_3))=s_3$. Let the set $W=\{S(s_2,\theta);
\theta\in[\lambda(r), \pi-\lambda(r)]\}$. We choose $\ve_1>0$ such
that $\nu(\Pi(r))- \ve_1\geq\ve_2>0$. Then for $\nu(A)\leq \ve_1$ we
have $\nu(\Pi(r)\backslash A)\geq \ve_2$. We first consider the case
when $\nu((W\times [r,\theta-r])\cap (\Pi(r)\backslash
A))\geq\ve_3>0$. This implies that there exist $\ve_4>0$ such that
$\text{mes}(\{s_3\in W; (s_3,\theta_3)\in \Pi(r)\backslash A \text{
for some } \theta_3\in [r, \pi-r]\})\geq\ve_4>0$. Here
$\text{mes}(\cdot)$ is the standard Lebesgue measure on
$\mathbb{R}$. We have, for some $\ve_5=\ve_5(\lambda(r))>0$
depending on $\lambda(r)$, that $\left|\dfrac{\partial
\zeta}{\partial s_3}\right|\geq \ve_5>0$ for $s_3\in W$. Let
$K=\{\zeta(s_2,s_3); s_3\in W, (s_3,\theta_3)\in \Pi(r)\backslash A
\text{ for some } \theta_3\in [r, \pi-r]\}$. We see that
$\text{mes}(K)\geq\ve_4\ve_5>0$. We also notice that $K \subseteq
[\lambda(r), \pi-\lambda(r)]$. Let us define the set
$M=\{\theta_3\in [r, \pi-r]; (s_3,\theta_3)\in \Pi(r)\backslash A
\text{ for some } s_3\in W\}$. There exist $\ve_6>0$ such that
$\text{mes}(M)\geq \ve_6>0$. Now we have

$$\begin{array}{l}
p^{(2)}(x,\Pi(r)\backslash A) \\
\play{=\int_{\theta_3\in M} \int_{\theta \in K} p_I(\delta;
\theta_2, \theta)p_I(\delta; \Theta(s_2,\theta), \theta_3)d\theta
d\theta_3}\\
\play{\geq \int_{\theta_3\in M} \int_{\theta \in K}
a(\lambda(r))a(\lambda(\lambda(r)))d\theta d\theta_3}
\\
\play{\geq a(\lambda(r))a(\lambda(\lambda(r)))\ve_4\ve_5\ve_6>0 \ .}
\end{array}$$

Therefore $p^{(2)}(x,A)\leq
1-a(\lambda(r))a(\lambda(\lambda(r)))\ve_4\ve_5\ve_6$. So in this
case we can choose $n_0=2$ and
$0<\ve<\min(a(\lambda(r))a(\lambda(\lambda(r)))\ve_4\ve_5\ve_6,\ve_1)
$.

Now let us consider the case when $\nu((W\times [r,\theta-r])\cap
(\Pi(r)\backslash A))=0$. Since we have $\nu(\Pi(r)\backslash A)\geq
\ve_2>0$, we have $\nu((([0,L]\backslash W)\times [r,\theta-r])\cap
(\Pi(r)\backslash A))\geq \ve_2>0$. In this case we can take one
more transition from $(s_2,\theta_2)$ to some other set of positive
measure lying away from $[0,L]\backslash W$ (provided that $r$ is
small enough) and carry out a similar argument as above. We conclude
in this case that $n_0=3$.

It is easy to see that this chain has only one ergodic component
$\Pi(r)$, since from any point $x\in \Pi(r)$ the chain can come to
any neighborhood of any other point $y \in \Pi(r)$ in several steps.

Since $m(x)$ is an invariant density for both the billiard map $f$
and the diffusion process $I_t^{\theta}$, it is also invariant for
the chain $\{X_n^\delta\}_{n \geq 1}$. Therefore, after
normalization, the unique invariant density for the chain
$\{X_n^\delta\}_{n \geq 1}$ is $d(x)=\play{\frac{1}{2L}m(x)}$.
$\blacksquare$

\

In particular, Law of Large Numbers holds for this Markov chain
$\{X_n^\delta\}_{n \geq 1}$ on $\Pi$.

\

\textit{Proof of Lemma 5.3:} Suppose there are $n=n(\varepsilon)$
collisions happened during time
$\displaystyle{[\frac{t}{\varepsilon}, \frac{t+\Delta
t}{\varepsilon})}$, and the corresponding successive positions of
the particle on $\Pi$ are $X_{k+1}^\delta,...,X_{k+n}^\delta$. It is
clear that $n \rightarrow \infty$ as $\varepsilon \downarrow 0$ and
$n\varepsilon \sim O(\Delta t)$ as $\varepsilon \downarrow 0$. For
small $\Delta t
> 0$, we have

$$H^{\ve,\delta}_{(t+\Delta
t)/\varepsilon}-H^{\ve,\delta}_{t/\varepsilon} =-\Delta t
(\frac{\Delta
t}{n\varepsilon})^{-1}\displaystyle{\frac{c(X_{k+1}^\delta)+c(X_{k+2}^\delta)+...+c(X_{k+n}^\delta)}{n}}.
$$

We also have $$\displaystyle{\frac{\Delta t}{n
\varepsilon}=\frac{t_{initial}-t_{final}+T(X_{k+1}^\delta,H_1)+T(X_{k+2}^\delta,H_2)
+...+T(X_{k+n}^\delta,H_n)}{n}}, \eqno(5.7)$$ where $H_i \ , \ 1
\leq i \leq n$ are the energies of the particle at the $i$-th
collision (before it loses some energy); $t_{initial}$ is the time
between the starting time $\displaystyle{\frac{t}{\varepsilon}}$ of
the motion and the time of the first collision with $\partial G$;
$t_{final}$ is the time between the ending time
$\displaystyle{\frac{t+\Delta t}{\varepsilon}}$ of the motion and
the time after that when the particle has next (the $n+1$-th)
collision with the boundary of $G$. Since function $T(x,H)$ is
Lipshitz in $H$ (when $H>h>0$) with respect to $x$, we have
$|T(x,H_i)-T(x, H^{\varepsilon}_{t/\varepsilon})|\leq
C|H_i-H^{\varepsilon}_{t/\varepsilon}|\sim O(\Delta t)$ as
$\varepsilon \downarrow 0$. Therefore from (5.2) , by using the Law
of Large Numbers for this Markov chain , we get
$$\frac{\Delta t}{n\varepsilon}\ra\mathbb{E}T(X,H^{\varepsilon}_{t/\varepsilon})+O(\Delta t)$$
in probability, as $\ve \da 0$.

By using the Law of Large Number for our Markov chain we also see
that
$$\frac{c(X_{k+1}^\delta)+c(X_{k+2}^\delta)+...+c(X_{k+n}^\delta)}{n}\ra\mathbb{E}c(X)$$
in probability, as $n \ra \infty$.

Here the random variable $X$ has the stationary distribution of this
chain on $\Pi$ with density $d(x)=\play{\frac{1}{2L}m(x)}$.
Therefore $$\lim_{\varepsilon \downarrow 0}\frac{H^{\varepsilon,
\delta}_{(t+\Delta t)/\varepsilon}-H^{\varepsilon,
\delta}_{t/\varepsilon}}{\Delta
t}=-\frac{\mathbb{E}c(X)}{\mathbb{E}T(X,H^{\varepsilon}_{t/\varepsilon})+O(\Delta
t)}$$ in probability.

This means that as $\varepsilon \downarrow 0$, the rescaled process
$H^{\ve,\delta}_{t/\varepsilon}$ within edge $I_3$ of the tree
$\Gamma$ converges in probability to the trajectory of the
deterministic process defined by the differential equation

$$\displaystyle{\frac{dH_t^\delta}{dt}=-\frac{\displaystyle{\iint_{\Pi}c(x)m(x)dx}}{\displaystyle{\iint_{\Pi}T(x,H_t)m(x)dx}}}.$$

Notice that the this is independent of $\delta$. Therefore letting
$\delta \downarrow 0$, we get the desired result. $\blacksquare$

\

\textit{Proof of Lemma 5.4:} This is a known result in integral
geometry (see [11]). But we will still provide a short proof here in
order for the reader to understand some calculations later in this
section.

\begin{figure}
\includegraphics[height=9cm, width=15cm , bb=43 11 481 273]{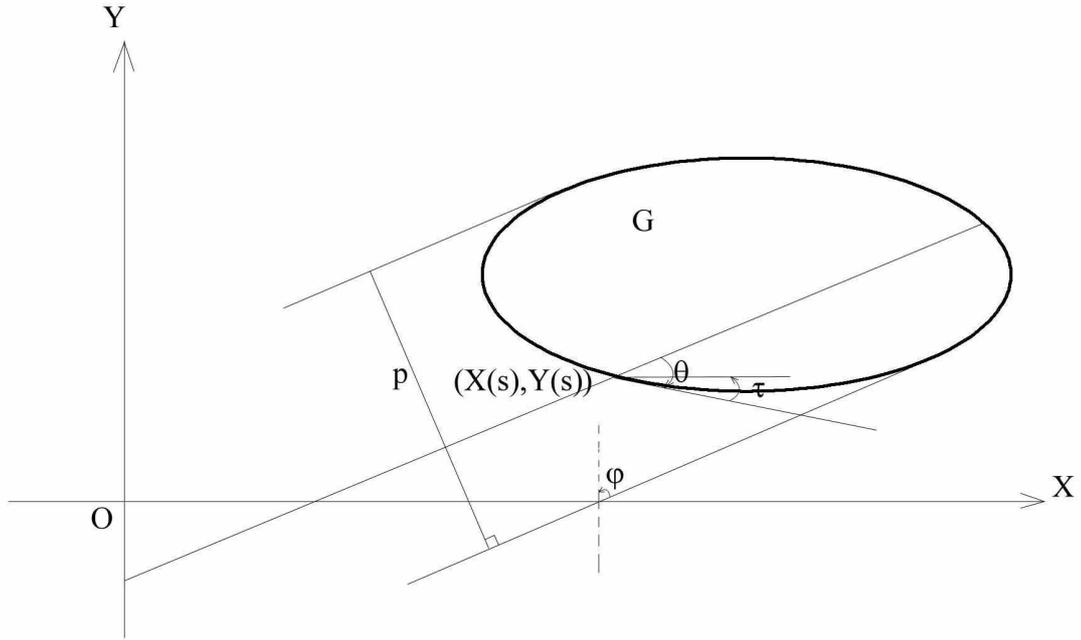}
\caption{A problem in integral geometry}
\end{figure}

Let the equation of the line intersecting the closed curve $\partial
G$ be written in the form (see Fig.12)

$$X\cos \varphi + Y \sin \varphi=p \ .$$

Let the curve $\partial G$ be $$X=X(s), Y=Y(s)$$ where $s$ is the
arc length parameter.

Let $\tau$ be the angle of the tangential direction at point
$(X,Y)$, we have $$\cos \tau =X'(s) \  , \ \sin \tau = Y'(s) \ , \
\frac{\pi}{2}+\varphi-\tau=\theta \ .$$

We then have

$$
\begin{array}{l}
d\varphi=d\theta+d\tau=d\theta+\kappa ds \ .
\end{array}$$

Here $\kappa=\kappa(X,Y)$ is the curvature at point $(X,Y)$, and

$$\begin{array}{l}
dp=X'(s)\cos\varphi ds - X(s) \sin \varphi d\varphi + Y'(s) \sin
\varphi ds + Y(s) \cos \varphi d\varphi\\ \ \  \ \ =(X'(s)\cos
\varphi + Y'(s)\sin \varphi) ds -X(s)\sin \varphi(\kappa
ds+d\theta)+Y(s)\cos\varphi (\kappa ds + d\theta)\\ \  \ \ \
=(X'(s)\cos \varphi + Y'(s)\sin \varphi+\kappa(-X(s)\sin
\varphi+Y(s)\cos\varphi))ds+(-X(s)\sin
\varphi+Y(s)\cos\varphi)d\theta \ .
\end{array}$$

These equalities imply that

$$\begin{array}{l}
dp\wedge d\varphi
\\=(X'(s)\cos \varphi + Y'(s)\sin
\varphi+\kappa(-X(s)\sin \varphi+Y(s)\cos\varphi))ds\wedge d\theta
\\ \ \ \ + \kappa(-X(s)\sin \varphi+Y(s)\cos\varphi)d\theta\wedge ds
\\=(X'(s)\cos \varphi+Y'(s)\sin \varphi)ds\wedge d\theta
\\=(\cos \tau \cos \varphi+\sin \tau \sin \varphi)ds\wedge d\theta
\\=\cos(\tau-\varphi)ds\wedge d\theta = \sin\theta ds \wedge d\theta \ .
\end{array}$$

So we have that $$\int_0^\pi\int_0^L L(s,\theta)\sin \theta ds\wedge
d\theta=2\int_0^\pi\int_0^P L(s,\theta)dp \wedge d\varphi=2\pi A.$$

We have used the fact that $\play{\int_0^P L(s,\theta)dp=A}$, where
$P$ is the height of $G$ projected onto the direction perpendicular
to $l(s,\theta)$ (recall that $l(s,\theta)$ is the chord starting
from $(s,\theta)\in \Pi$). $\blacksquare$

\

We would like to mention that the relation $dp\wedge d\varphi=\sin
\theta ds \wedge d\theta$ will be used in some of our later
calculations. In future calculations, we will use the abbreviated
form $dpd\varphi=\sin\theta dsd\theta$.

\

\textit{Proof of Lemma 5.5:} From Lemma 5.2 we already know that the
embedded chain $\{X_n^\delta\}_{n \geq 1}$ has a unique invariant
measure on $\Pi$. Since our process $\Lambda_t^{\ve,\delta}$ is
purely deterministic within $\bigwedge \setminus \Pi$, and any point
in $\bigwedge \setminus \Pi$ can be reached in finite time from a
point on $\Pi$, $\Lambda_t^{\ve,\delta}$ must be positive recurrent
on $\bigwedge$ and therefore it has a unique invariant measure on
$\bigwedge$.

We use the Khasminskii formula (see, [6], on page 185, formula
(4.1)) for the expression (up to a constant factor) of the invariant
measure $\widehat{\mu}(\cdot)$ for $\Lambda_t^{\ve,\delta}$ on
$\bigwedge$ through the invariant measure $\dfrac{1}{2L}m(x)dx$ for
the embedded chain $\{X_n^\delta\}_{n \geq 1}$ on $\Pi$: for any
measurable set $B\subseteq \bigwedge$,

$$\widehat{\mu}(B)=\frac{1}{2L}\iint_\Pi m(x)dx \mathbb{E}_x\int_0^{\tau_1^x}\chi_B(\Lambda_t^{\ve,\delta})dt \ . \eqno(5.8)$$

Here $\tau_1^x$ is the first time when $\Lambda_t^{\ve,\delta}$,
$\Lambda_0^{\ve,\delta}=x\in \Pi$, reaches $\Pi$ again. Note that
for a given $x$, $\tau_1^x$ as well as the motion
$\Lambda_t^{\ve,\delta}$, $0<t<\tau_1^x$, and in particular, the
time spent by $\Lambda_t^{\ve,\delta}$ in $B$, are not random, so
that the expectation sign in (5.8) can be omitted. Since the
invariant measure of our imbedded chain and of the billiard map are
the same, we conclude from (5.8) that $\widehat{\mu}(B)$ coincides
with the invariant measure of the non-perturbed billiard. As is
known, the Lebesgue measure is invariant for the billiard. Thus
$\widehat{\mu}(B)$ is the uniform distribution on $\bigwedge$.
$\blacksquare$

\

\textit{Proof of Lemma 5.6:} This follows from the fact that

$$P\{\max\li_{0\leq t <T<
\infty}|H_t^{\ve,\delta}-\widehat{H}_t^{\ve,\delta}|<C\ve\}=1$$

where $C>0$ is a constant. $\blacksquare$

\

\textit{Proof of Lemma 5.7:} Consider a time change in the process
$\widehat{N}_t^{\ve,\delta}$. If at time $t$,
$\Lambda_t^{\ve,\delta}$ lies on a chord $l(x)$ starting from
$x=(s,\theta)\in \Pi$, then we set
$d\widetilde{t}=\play{\frac{L(x)}{c(f(x))}}dt$. Here
$f(x)=f(s,\theta)=(S(s,\theta), \Theta(s,\theta))$ is the billiard
map and $L(x)$ is the length of the chord $l(x)$. The density of the
invariant measure for
 $\Lambda_{\widetilde{t}}^{\ve,\delta}$ (with respect to Lebesgue measure) is proportional to
 $\play{\frac{c(f(x))}{L(x)}}$. Here $x=x(g,\varphi)$ is a point in $\Pi$ such
 that $(g,\varphi)\in \Lambda$
 lies on the chord $l(x)$ starting from $x$. This time change
 corresponds to multiplication of the generator of the process
 $\widehat{N}_t^{\ve,\delta}$ by $\play{\frac{L(x)}{c(f(x))}}$.

The $H$-component of the process obtained after time change moves
down uniformly with a non-random speed, so that it hits level $H(O)$
at a non-random time $t=t(H_0)$. Let $\Lambda_1=\{(g,\varphi)\in
\bigwedge : f(x(g,\varphi))\in \Pi_1\}$, $\Lambda_2=\{(g,\varphi)\in
\bigwedge : f(x(g,\varphi))\in \Pi_2\}$. Here $\Pi_1, \Pi_2$ are
defined as before, i.e., $\Pi_i, \ i \in \{1,2\}$ is the set of
those points on $\Pi$ which correspond to well $i \ , \ i \in
\{1,2\}$.

By using Lemma 5.5, we see that

\

$\lim\li_{\ve \da 0} P\{ \widehat{N}_{\widetilde{t}}^{\ve,\delta}$
finally falls into well 1 $\}$
$$\begin{array}{l}
=D\play{\iiint_{\Lambda_1}\frac{c(f(x(g,\varphi)))}{L(x(g,\varphi))}dgd\varphi}
\\=D\play{\iiint_{\Lambda_1}\frac{c(f(x(g,\varphi)))}{L(x(g,\varphi))}dL(x(g,\varphi))dpd\varphi}
\\=D\play{\iint_{f^{-1}(\Pi_1)}\frac{c(f(s,\theta))}{L(s,\theta)}\left(\int_0^{L(s,\theta)}dL(s,\theta)\right)\sin{\theta}dsd\theta}
\\=D\play{\iint_{f^{-1}(\Pi_1)}c(f(s,\theta))m(s,\theta)dsd\theta}
\\=D\play{\iint_{\Pi_1}c(s,\theta)m(s,\theta)dsd\theta}.
\end{array}$$

Here $m(s,\theta)=\sin \theta$ is the density of the Liouville
measure. We have used the fact that
$dgd\varphi=dL(s,\theta)dpd\varphi=dL(s,\theta)\sin\theta dsd\theta$
(notice here that $dg=dLdp$), as was calculated in Lemma 5.4, and
the invariance of Liouvuille measure under the billiard map.
Normalizing constant
$D=\left(\play{\iiint_{\bigwedge}\frac{c(f(x(g,\varphi)))}{L(x(g,\varphi))}dgd\varphi}\right)^{-1}
=\left(\play{\iint_{\Pi}c(s,\theta)m(s,\theta)dsd\theta}\right)^{-1}$.
Similarly, we have $\lim\li_{\ve \da 0} P\{
\widehat{N}_{\widetilde{t}}^{\ve,\delta}$ finally falls into well 2
$\}=\play{\frac{\play{\iint_{\Pi_2}c(s,\theta)m(s,\theta)dsd\theta}}{\play{\iint_{\Pi}c(s,\theta)m(s,\theta)dsd\theta}}}.
$

The trajectory of $\widehat{N}_{\widetilde{t}}^{\ve,\delta}$ is the
same as $\widehat{N}_t^{\ve,\delta}$. Taking into account this fact
and Lemma 5.6, the statement follows. $\blacksquare$

\

\

\textbf{Acknowledgement:} This work was supported in part by NSF
grants DMS-0803287 and DMS-0854982. We would like to thank the
anonymous referee for carefully reading the manuscript.

\

\

\section*{References}

[1] V.I.Arnold, \textit{Mathematical Methods of Classical
Mechanics}, Springer, 1978.

[2] M.Brin, M.Freidlin, On stochastic behavior of perturbed
Hamiltonian systems, \textit{Ergodic Theory and Dynamical Systems},
1999/2000, \textbf{20}: 55-76.

[3] J.L.Doob, \textit{Stochastic Processes}, New York, John Wiley,
1953.

[4] M.I.Freidlin, \textit{Functional Integration and Partial
Differential Equations}, Princeton University Press, 1985.

[5] M.I.Freidlin, M.Weber, Random perturbations of dynamical systems
and diffusion processes with many conservation laws, \textit{Probab.
Theory Relat. Fields}, \textbf{128},441-466, 2004.

[6] M.I.Freidlin, A.D.Wentzell, \textit{Random perturbations of
dynamical systems}, Springer, 1998.

[7] I.I.Gikhman, A.V.Skorohod, \textit{The theory of stochastic
processes}, Vol I, Springer, 2004.

[8] K.It\^{o}, H.P.MacKean, \textit{Diffusion processes and their
sample paths}, Springer, 1974.

[9] A.Katok, B.Hasselblatt, \textit{Introduction to the modern
theory of dynamical systems}, Encyclopedia of mathematics and its
applications, Vol 54, Cambridge University Press, 1995.

[10] V.V.Petrov, \textit{Sums of Independent Random Variables},
Springer, 1975.

[11] L.A.Santalo, \textit{Introduction to Integral Geometry}, Paris,
1953.

\end{document}